\documentclass[final]{siamltex}
\pdfoutput=1

\usepackage{graphicx,color,enumerate,nicefrac,bm}
\usepackage{algorithmicx,algpseudocode}

\usepackage[T1]{fontenc}
\usepackage[utf8]{inputenc}
\usepackage[dvipsnames,svgnames,table]{xcolor}
\usepackage{listings}
\usepackage[scaled=0.8]{beramono}

\definecolor{LightGrey}{rgb}{0.9629411,0.9629411,0.9629411}
\definecolor{LighterGrey}{gray}{0.99}
\definecolor{Mauve}{rgb}{0.58,0,0.82}

\usepackage{graphicx}
\usepackage{tikz}
\usepackage{moredefs}
\usepackage{float}
\usepackage{url}
\usepackage{mathtools}
\usepackage{amsfonts}
\usepackage{xfrac}
\usepackage{hhline}

\usepackage[super]{nth}
\usepackage{caption}
\usepackage{subcaption}
\usepackage{algpseudocode, algorithm}

\usepackage{booktabs,siunitx}

\usepackage{wrapfig}
\usepackage{pdflscape}
\usepackage{rotating}
\usepackage{chngpage}

\usepackage{pgfplots}
\usepackage{pgfplotstable}
\usepackage{longtable}
\usepackage{array}
\usepackage{booktabs}

\usepackage{ifthen}

\pgfplotsset{compat=newest} 
\pgfplotsset{plot coordinates/math parser=false} 
\newlength\figureheight 
\newlength\figurewidth 
\usetikzlibrary{3d,calc}
\usetikzlibrary{intersections}
\usepackage{tkz-euclide}
\usetkzobj{all}

\usepackage{tikz-qtree}
\usetikzlibrary{trees}
\usetikzlibrary{%
    calc,%
    fadings,%
    shadings%
}
\usetikzlibrary{arrows,snakes,shapes}
\usetikzlibrary{backgrounds}
\usetikzlibrary{
    shapes.geometric,
    decorations.pathreplacing
}
\usetikzlibrary{positioning,shadows,arrows}

\tikzstyle{block} = [rectangle, draw, text width=1cm, text centered, rounded corners, minimum height=4em]

\newcommand{\nno}{\nonumber}
\def\hh{{\tt h}}
\newcommand{\Cip}{C_{\mbox{\tiny IP}}}

\newcommand{\cvec}[1]{\mathbf{#1}}
\newcommand{\average}[1]{\{\!\!\{{#1}\}\!\!\}}

\newcommand{\tenjump}[1]{\ensuremath{\underline{[\![#1]\!]}}}

\newcommand{\dd}[1]{\ensuremath{\;\mathrm{d} #1}}

\newcommand{\defeq}{\vcentcolon=}

\newcommand{\uu}[1]{\mathbf{#1}}
\newcommand{\uuu}[1]{\underline{#1}} 
\newcommand{\mesh}{{\mathcal T}_h}

\newcommand{\ilc}[1]{\lstinline[language=Python]{#1}}

\newcommand{\visc}[1][]{\ensuremath{\mu_{#1}}}

\newcommand{\prandtl}{\ensuremath{\mathrm{Pr}}}

\newcommand{\omegatri}{\ensuremath{\mathcal{T}_{h}}}

\newcommand{\domain}[1][]{\ensuremath{\Omega}_\text{#1}}

\newcommand{\eyedentity}{\ensuremath{\underline{\mathbf{I}}}}

\newcommand{\extbdr}[1][]{{\Gamma_{#1}}}
\newcommand{\intbdr}{\ensuremath{{\Gamma_\mathcal{I}}}}

\newcommand{\vvec}[2]{\ensuremath{\mathbf{V}^{#2}_{\polyo}(#1)}}

\newcommand{\polyo}{\ensuremath{\ell}}

\DeclareSIUnit\torr{torr}

\usepackage{color}

\title{Automatic Symbolic Computation for
Discontinuous Galerkin Finite Element Methods}

\author{
Paul Houston
\thanks{
School of Mathematical Sciences, University of Nottingham,
University Park, Nottingham, NG7 2RD, UK,
email: {\tt Paul.Houston@nottingham.ac.uk}.}
\and Nathan Sime
\thanks{Department of Engineering, University of Cambridge,
Trumpington Street,
Cambridge CB2 1PZ, UK,
email: {\tt njcs4@cam.ac.uk}.}
}

\begin{document}

\maketitle

\begin{abstract}
The implementation of discontinuous Galerkin finite element methods  (DGFEMs)
represents a very challenging computational task, particularly for
systems of coupled nonlinear PDEs, including multiphysics problems, 
whose parameters may consist of power series or functionals of the 
solution variables. Thereby, the exploitation of symbolic algebra to 
express a given DGFEM approximation of a PDE problem
within a high level language, whose syntax
closely resembles the mathematical definition, is
an invaluable tool. Indeed, this then facilitates the automatic assembly
of the resulting system of (nonlinear) equations, as well as the computation
of Fr\'{e}chet derivative(s) of the DGFEM scheme, needed, for example, within
a Newton-type solver. 
However, even exploiting symbolic algebra, the discretisation of coupled systems
of PDEs can still be extremely verbose and hard to debug. Thereby,
in this article we develop a further layer of abstraction by designing
a class structure for the automatic computation of DGFEM formulations. 
This work has been implemented within the
FEniCS package, based on exploiting the Unified Form Language. Numerical
examples are presented which highlight the simplicity of implementation of
DGFEMs for the numerical approximation 
of a range of PDE problems. 
\end{abstract}



\pagestyle{myheadings} \thispagestyle{plain} \markboth{P. HOUSTON and N. SIME}{Symbolic Computation for DGFEMs}

\section{Introduction}

The finite element method (FEM) represents an indispensable computational tool
for the accurate, efficient, and rigorous numerical approximation of continuum
models arising within a wide range of scientific and engineering application
areas. Key reasons for the success of FEMs include their applicability to very
general classes of partial differential equations (PDEs), simple treatment of
complicated computational geometries and enforcement of boundary conditions,
ease of adaptivity including both local mesh subdivision ($h$--refinement) and
local polynomial enrichment ($p$--refinement), and, from a mathematical point of
view, the availability of tools for their rigorous error analysis. However, when
compared with their finite difference counterparts, FEMs are typically regarded
as being {\em complicated} to implement. Indeed, assembly of the underlying
matrix stemming from the FEM discretisation of a given (linear, for example) PDE
problem typically involves mapping each element present in the computational
mesh, defined in the global coordinate system, to a given reference or canonical
element, where both the local FEM basis and corresponding quadrature, needed to
approximate the underlying integral, is defined. In this manner, local elemental
stiffness matrices and load vectors may be computed; these entries are then
inserted into the global matrix and right-hand side vector, respectively,
according to the elementwise local-to-global degree of freedom mapping, subject
to the enforcement of inter-element continuity constraints and the imposition of
boundary conditions. This general assembly strategy is almost universally
employed within both open source and commercial software, e.g., 
FreeFem++~\cite{freefem}, DUNE~\cite{ans26526}, deal.II~\cite{BangerthHartmannKanschat2007}
and OpenFOAM~\cite{Weller1998}, to name just a few. 
Thereby, in principle, assuming that a user can define both the
element stiffness matrix and load vector, then a numerical approximation may be
readily determined. However, most open source software packages do not provide
easy-to-use user interfaces; indeed, typically the user must first understand
the low level language in which the code is written, e.g., C++, Fortran, etc,
and develop an understanding of the underlying datastructures and
function/subroutine calls defined within their chosen package in order to be
able to develop code specific to their own application. This can, of course, be
a rather time-consuming exercise, and often requires one to use multiple
software libraries, for example, when certain features are not available within
the chosen package which the user needs to utilise.

In the nonlinear setting, assuming a Newton-type iteration is exploited to solve
the underlying system of nonlinear equations stemming from the given FEM
employed, the general strategy of assembly of the resulting linearised equations
is similar to the linear case; in this setting the load vector represents the
residual of the numerical scheme. However, in this case, the Fr\'{e}chet
derivative of the FEM must now be computed; we point out that, in the context of
the numerical detection of bifurcation points, the number of derivatives that
must be computed just to form the so-called {\em extended system} needed to
accurately compute the bifurcation point, even before a Newton iteration is
implemented, is dependent on the codimension of the singularity being sought,
cf.~\cite{Cliffe2000a,Cliffe2000}. In general, the evaluation and
implementation of both the FEM residual and, moreover, the Fr\'{e}chet
derivative(s) of the scheme is a difficult and time consuming task, which is
inevitably prone to human error. This is particularly the case when exploiting
FEMs for the numerical approximation of systems of coupled nonlinear PDEs,
including multiphysics problems, whose parameters may consist of power series or
functionals of the solution variables, cf.~\cite{cvd_paper}.

Thus far, our discussion has primarily focussed on the application of conforming
FEMs, whereby, excluding Neumann, or weakly imposed boundary
conditions, only element contributions need to be evaluated. The exploitation of
more general FEMs, and in particular discontinuous Galerkin FEMs (DGFEMs)
requires the implementation of inter-element flux terms, which involves
combinations of both inner-- and outer--traces of the FEM solution, relative to
a given orientation of the element face. In recent years, DGFEMs have become an
increasingly popular class of FEMs, most notably due to their local conservation
properties, inherent numerical stability for convection--dominated diffusion
problems, limited interelement communication, which is restricted only to
neighbouring elements, and has important advantages for the implementation of
boundary conditions and the parallel efficiency of the method, and finally
the ease in which so--called hanging nodes can be treated, and the efficient
implementation of $hp$--adaptivity. Indeed, tremendous progress has been made on
both the analytical and computational aspects of DGFEMs; for a review of some of
the key developments in the subject, we refer to the recent
monographs~\cite{dg_book,DiPietroErn,hesthaven_DG,riviere_DG}. For a historical
review of DGFEMs, we refer to the articles~\cite{unified,MR1728854}, and the
references cited therein.

The addition of inter-element face terms within DGFEMs further complicates the
implementation of such schemes within standard FEM packages. Assuming for a
moment that the underlying PDE problem is linear, we note that for a given
interior face, shared by two neighbouring elements, there are, in general, four
local face matrices, which stem from the different combinations of traces of
the FEM solution from the two elements whose boundaries form the current face.
Once these local face--wise matrices have been assembled, they can then be
inserted into the global FEM matrix in an analogous manner to the treatment of
the element stiffness matrix. On boundary faces, contributions to the load
vector must also be computed. Again, in the nonlinear setting, the task of
computing both the residual vector and Fr\'{e}chet derivative(s) 
of the DGFEM scheme is a very technical and time consuming task.

The purpose of this article is to discuss the use of symbolic algebra to
facilitate the assembly of FEM matrix problems, and in particular those arising
from the application of DGFEMs, for the numerical approximation of general
nonlinear systems of PDEs. The general approach is to develop a high level
language syntax, which closely corresponds to the mathematical formulation of
the underlying FEM. Thereby, through this layer of abstraction, the user needs
to only specify the FEM residual in a concise and easy to read manner, whereby
the evaluation of the Fr\'{e}chet derivative(s) of the scheme are automatically
computed symbolically and the resulting low level C++/Fortran code for element
stiffness and face matrices and load/residual vectors are automatically generated
by the so-called {\em form compiler}. In
this way, any existing open source FEM package may be utilised, subject to the
implementation of a suitable interface which directly calls the automatically
generated snippets of code when assembling element and face matrices and load/residual
vectors. Most importantly, we stress that once the user has selected a
particular FEM package which is appropriate for their purposes, the interface to
that software platform which links to the automatically generated code only
needs to be written and debugged {\em once}; it may then subsequently be
exploited to solve a potentially huge variety of PDE problems, with a plethora
of FEM schemes.
At this point it is pertinent to mention that some FEM packages do indeed
include such a symbolic interface at the heart of their design; most notably we
mention the excellent FEniCS package, cf.~\cite{AlnaesBlechta2015a,Logg2012},
for example. The Unified Form Language (UFL) component of FEniCS provides an
easy to use python interface which allows for the automatic FEM numerical
approximation of systems of PDEs in a user-friendly manner. Other
form compilers include SyFi~\cite{Simula.simula.1092}, and
Manycore Form Compiler~\cite{MARKALL20101815}, for example.

However, even exploiting a package such as UFL, as
powerful as it is, the definition of DGFEMs in this framework is still rather
verbose, particularly for nonlinear systems of coupled PDEs. With this in mind,
we present a further layer of abstraction for the automatic computation of DGFEM
formulations employing symbolic algebra. Indeed, in this article DGFEM utility
functions for general PDE operators are developed, which significantly simplifies
the specification of the DGFEM discretisation of a given problem. This
work was originally inspired by the need to numerically model the formation of a
hydrogen plasma in a microwave power assisted chemical vapour
deposition reactor employed for the manufacture of synthetic
diamond; this includes constituent equations for   
the background gas mass average velocity, gas
temperature, electromagnetic field energy and plasma density, 
cf.~\cite{cvd_paper,nate_phd}.



In \cite{nate_phd} we originally developed easy-to-use DGFEM utility functions
as part of our inhouse software package AptoPy~\cite{nate_phd}, for application
with our own FEM package AptoFEM~\cite{aptofem}. AptoPy is written in Python and
exploits the open source symbolic algebra package SymPy~\cite{sympy}. We stress
that these choices are entirely user-dependent; indeed, in the past we have
employed the symbolic algebra packages REDUCE~\cite{reduce}, Mathematica, and
Maple. In addition to AptoFEM, ENTWIFE has also been employed as the low level
FEM package. The DGFEM utility functions written in AptoPy have been ported to
UFL for use within the FEniCS package; full open source codes are available from
\url{https://bitbucket.org/nate-sime/dolfin_dg}. With this in mind, for
simplicity of presentation, throughout this article, we only show snippets of
UFL code in order to highlight how the DGFEM utility functions may be exploited
in practice; the corresponding AptoPy syntax is quite similar,
cf.~\cite{nate_phd}.

The outline of this article is as follows. In Section~\ref{sec:dg-fems}
we briefly outline the DGFEM discretisation of general systems of 
nonlinear conservation laws. Then, in Section~\ref{sec:comp_framework_for_dgfems}
we propose a computational framework for the automatic generation of
DGFEM schemes within a simple unified setting. On the basis of this work,
in Section~\ref{sec:examples} we provide some examples to illustrate the
flexibility and ease within which systems of PDEs may be numerically 
approximated using the software developed in this article. 
Finally, in Section~\ref{sec:con_rem}
we provide a summary of the work undertaken in this article, as well as
outlining potential future developments.

\section{Discontinuous Galerkin Finite Element Methods}
\label{sec:dg-fems}

As a representative PDE example, in this section we outline the DGFEM
discretisation for the following system of conservation laws:
\begin{equation}
\nabla\cdot \left({\mathcal F}^c({\bf u})
-{\mathcal F}^v({\bf u}, \nabla{\bf u})\right) = {\bf f} \quad\mbox{in } \Omega,
\label{eq:model-prob}
\end{equation}
where $\Omega$ is an open bounded domain in ${\mathbb R}^d$, $d \geq 1$, with
boundary $\Gamma$. Here, ${\bf u} = (u_1,\ldots,u_m)^\top$, $m \geq 1$, 
${\mathcal F}^c({\bf u}) = (\mathbf{f}^c_1(\cvec{u}), \ldots,
\mathbf{f}^c_d(\cvec{u})) $ and ${\mathcal F}^v({\bf u}, \nabla{\bf u}) =
(\cvec{f}^v_1(\cvec{u}, \nabla \cvec{u}), \ldots,$ $\cvec{f}^v_d(\cvec{u}, \nabla
\cvec{u}) )$ represent the convective and diffusive fluxes, respectively, which
are assumed to be continuously differentiable, and ${\bf f}$ is a given
source function. For simplicity of presentation,
we assume that~\eqref{eq:model-prob} may be supplemented with the boundary
conditions:
\begin{eqnarray}
		\mathbf{u}          & = &\mathbf{g}_D \; ~~~ \text{on} \; ~ \extbdr_D, \label{eq:bcs1}\\
	\mathcal{F}^v(\mathbf{u}, \nabla \mathbf{u}) \cdot \mathbf{n} &=& \mathbf{g}_N
	\; ~~~ \text{on} \; ~\extbdr_N,
	\label{eq:bcs2}
\end{eqnarray}
where $\Gamma=\extbdr_D\cup \extbdr_N$, and $\extbdr_D$ and $\extbdr_N$ are two
disjoint subsets, with $\extbdr_D$ nonempty and relatively open in $\Gamma$. We
stress that more general boundary conditions can also be considered; for
example, in the case when~\eqref{eq:model-prob} represents the compressible
Euler or Navier-Stokes equations, we refer to, for
example,~\cite{HartmannHouston01} and~\cite{HH08a}, respectively.

For the purposes of discretisation, we rewrite~\eqref{eq:model-prob} in the
following equivalent form:
\begin{equation*}
\nabla\cdot \left({\mathcal F}^c({\bf u})-G({\bf u})\nabla{\bf u}\right)\equiv\frac{\partial}{\partial x_k} \left( {\bf f}^c_k ({\bf u})-
G_{kl}({\bf u})\frac{\partial {\bf u}}{\partial x_l} \right) 
= {\bf f} \quad\mbox{in }   \Omega.
\end{equation*}
Here, the matrices 
\begin{equation}
G_{kl}({\bf u})=\partial{\bf f}^v_k({\bf u}, \nabla{\bf u})
/\partial u_{x_l}, \quad k,l=1,\ldots,d,
\label{eq:homo_tensor}
\end{equation}
are the homogeneity tensors defined by ${\bf f}^v_k({\bf u}, \nabla{\bf
u})=G_{kl}({\bf u})\partial{\bf u}/\partial x_l$, $k=1,\ldots,d$. We write the
homogeneity tensor product
\begin{equation}
    (G(\mathbf{u}) \nabla \mathbf{u})_{ik} = \sum_{j=1}^{m} \sum_{l=1}^{d}(G_{kl}\left(\mathbf{u}\right))_{ij} \frac{\partial \mathbf{u}_j}{\partial x_l}
    \label{eq:homo_tens_prod}
\end{equation}
such that $\mathcal{F}^v(\mathbf{u},\nabla \mathbf{u}) = G(\mathbf{u})\nabla \mathbf{u}$.

For simplicity of presentation, we now proceed to 
discretise~\eqref{eq:model-prob}--\eqref{eq:bcs2} based on employing the symmetric
interior penalty (SIPG) formulation presented in~\cite{HH08a}; however, we stress that other 
DGFEMs can easily be included within this general setting, cf. below. To this end, we partition
$\Omega$ into a mesh ${\mathcal T}_h=\{\kappa\}$ consisting of non-overlapping
open element domains~$\kappa$, such that 
$\bar{\Omega}=\cup_{\kappa\in {\mathcal T}_h}\bar{\kappa}$. 
For each $\kappa\in{\mathcal T}_h$, we denote by
$\uu{n}_{\kappa}$ the unit outward normal vector to the boundary $\partial
\kappa$. We assume that each $\kappa \in {\mathcal T}_h$ is an image of a fixed
reference element $\hat\kappa$, that is, $\kappa = \sigma_\kappa(\hat \kappa)$
for all $\kappa \in {\mathcal T}_h$. On the reference element $\hat\kappa$ we
define spaces of polynomials of degree $\polyo\geq 0$ as follows:
\begin{eqnarray}
{\mathcal Q}_{\polyo} = \mbox{span} 
 \left\{ \hat{{\bf x}}^{\alpha}: 
0 \leq \alpha_i \leq \polyo, ~0\leq i \leq d \right\}, \qquad
{\mathcal P}_{\polyo} = \mbox{span} 
      \left\{ \hat{{\bf x}}^{\alpha}: 0 \leq |\alpha| \leq \polyo \right\}.
\nno
\end{eqnarray}
Thereby, with this notation, we define the following DGFEM finite element space
\begin{eqnarray}
\vvec{\omegatri}{m}
= \{ {\bf v} &\in& \left[L_2(\Omega)\right]^m : 
  {\bf v} |_{\kappa} \circ \sigma_{\kappa} \in  \left[{\mathcal Q}_{\polyo}\right]^m
    \mbox{ if $\hat{\kappa} = \sigma_{\kappa}^{-1}(\kappa)$ is the unit hypercube, } \nno \\
&& \mbox{and } {\bf v} |_{\kappa} \circ \sigma_{\kappa} \in 
   \left[{\mathcal P}_{\polyo}\right]^m
    \mbox{ if $\hat{\kappa} = \sigma_{\kappa}^{-1}(\kappa)$ is the unit simplex};
    ~ \kappa \in {\mathcal T}_h \}. \nno
\end{eqnarray}

An {\em interior face} of ${\mathcal T}_h$ is defined as the (non-empty)
$(d-1)$--dimensional interior of $\partial \kappa^+ \cap\partial \kappa^-$,
where $\kappa^+$ and $\kappa^-$ are two adjacent elements of ${\mathcal T}_h$,
not necessarily matching. A {\em boundary face} $f$ of ${\mathcal T}_h$ is defined
as a (non-empty) $(d-1)$--dimensional facet of $\kappa$, $\kappa \in {\mathcal T}_h$,
where $\kappa$ is a boundary element of ${\mathcal T}_h$, such that 
$f\subset \partial \kappa \cap\Gamma$. We denote
by $\Gamma_{\mathcal I}$ the union of all interior faces of ${\mathcal T}_h$.
Let $\kappa^+$ and $\kappa^-$ be two adjacent elements of~ ${\mathcal T}_h$, and
$\uu x$ an arbitrary point on the interior face $f=\partial \kappa^+\cap\partial
\kappa^-$.  Furthermore, let $\uu{v}$ and~$\uuu{\tau}$ be vector- and 
matrix-valued functions, respectively, that are smooth inside each
element~$\kappa^\pm$. By $(\uu{v}^\pm,\uuu{\tau}^\pm)$, we denote the traces of
$(\uu{v},\uuu{\tau})$ on $f$ taken from within the interior of $\kappa^\pm$,
respectively.
Then, the averages of $\uu{v}$ and $\uuu{\tau}$ at $\uu{x}\in f$ are given by
$\average{{\bf v}} =({\bf v}^++{\bf v}^-)/2$ and $\average{\uuu{\tau}}
=(\uuu{\tau}^++\uuu{\tau}^-)/2$, respectively. Similarly, the jump of $\uu{v}$
at $\uu{x}\in f$ is given by $\tenjump{{\bf v}} ={\bf v}^+\otimes{\bf
n}_{\kappa^+}+{\bf v}^-\otimes{\bf n}_{\kappa^-}$, where we denote by
$\uu{n}_{\kappa^\pm}$ the unit outward normal vector of $\kappa^\pm$,
respectively. On $f \subset \Gamma$, we set $\average{{\bf v}} ={\bf v}$,
$\average{\uuu{\tau}} =\uuu{\tau}$ and $\tenjump{{\bf v}} ={\bf v}\otimes{\bf
n}$, where ${\bf n}$ denotes the unit outward normal vector to $\Gamma$.

The interior penalty DGFEM discretisation of~\eqref{eq:model-prob}--\eqref{eq:bcs2} is given by: find ${\bf u}_h \in \vvec{\omegatri}{m}$
such that
\begin{eqnarray}
\label{eq:dg-ns}
&& {\mathcal N}({\bf u}_h,{\bf v}_h) \nonumber \\
&&\equiv 
- \int_{\Omega} {\bf f}\cdot {\bf v}_h \dd \uu{x}
- \int_{\Omega} {\mathcal F}^c({\bf u}_h) : \nabla_h {\bf v}_h \dd \uu{x}
+ \sum_{\kappa \in {\mathcal T}_h} \int_{\partial \kappa\setminus\Gamma}
  {\mathcal H}({\bf u}_h^+,{\bf u}_h^-,{\bf n}^+) \cdot \uu{v}_h^+ \dd s \nonumber \\
&& 
\quad + \int_{\Omega} {\mathcal F}^v({\bf u}_h,\nabla_h{\bf u}_h) : \nabla_h{\bf v}_h \dd \uu{x}
- \int_{{\Gamma}_{\mathcal I}}  \average{{\mathcal F}^v({\bf u}_h,\nabla_h{\bf u}_h)}
  :\tenjump{{{\bf v}_h}} \dd s \nonumber \\
&& 
\quad - \int_{{\Gamma}_{\mathcal I}}
  \average{G^\top({\bf u}_h) \nabla_h{\bf v}_h} : \tenjump{{\bf u}_h} \dd s
 + \int_{{\Gamma}_{\mathcal I}}  \uuu{\delta}({\bf u}_h):\tenjump{{{\bf v}_h}} \dd s 
+{\mathcal N}_\Gamma({\bf u}_h, {\bf v}_h)  =0
\end{eqnarray}
for all ${\bf v}_h$ in $\vvec{\omegatri}{m}$. The subscript $h$ on the operator
$\nabla_h$ is used to denote the discrete counterpart of $\nabla$, defined
elementwise. Here, ${\mathcal H}(\cdot,\cdot,\cdot)$ denotes the (convective)
numerical flux function; this may be chosen to be any two--point monotone
Lipschitz function which is both consistent and conservative; typical choices
include the (local) Lax--Friedrichs flux, the HLLE flux, the Roe flux, 
and the Vijayasundaram flux,
cf.~\cite{Kro97,Tor97}.

The penalisation 
term $\uuu{\delta}(\cdot)$  arising in the DGFEM
scheme~(\ref{eq:dg-ns}) is given by
\begin{equation*}
  \uuu{\delta}({\bf u}_h) = ~\Cip \tfrac{\polyo^2}{\hh}\average{G({\bf u}_h)}\tenjump{{\bf u}_h}, 
\end{equation*}
where $\Cip$ is a (sufficiently large) positive constant, cf.~\cite{ghh_hpfem_part_II}. 
Moreover, $\hh \in L_\infty(\Gamma_{\mathcal I} \cup \Gamma)$ is
defined as $\hh({\bf x})=\min\{m_{\kappa^+},m_{\kappa^-}\}/m_f$, if ${\bf x}$ is
in the interior of $f=\partial \kappa^+\cap\partial \kappa^-$ for two
neighbouring elements in the mesh ${\mathcal T}_h$, and $\hh({\bf
x})=m_{\kappa}/m_f$, if ${\bf x}$ is in the interior of
$f=\partial\kappa\cap\Gamma$. Here, for a given (open) bounded set $\omega
\subset {\mathbb R}^s$, $s \geq 1$, we write $m_{\omega}$ to denote the
$s$--dimensional measure of $\omega$.

Finally, we define the boundary terms present in the form ${\mathcal
N}_\Gamma(\cdot,\cdot)$ by
\begin{eqnarray}
\label{eq:dg-ns-exterior}
{\mathcal N}_\Gamma({\bf u}_h, {\bf v}_h) &=& 
\int_\Gamma{\mathcal H}_\Gamma({\bf u}_h^+,{\bf u}_\Gamma({\bf u}_h^+),{\bf n}) \cdot {\bf v}_h^+ \dd s
    + \int_{\extbdr_D}\uuu{\delta}_\Gamma({\bf u}_h^+):{\bf v}_h \otimes {\bf n} \dd s \nno \\
&&    - \int_{\extbdr_N} \mathbf{g}_N \cdot {\bf v}_h \dd s
   -  \int_{\extbdr_D} 
     {\bf n}\cdot {\mathcal F}_{\Gamma}^{v}
     ({\bf u}_h^+,\nabla_h{\bf u}^+_h)\, {\bf v}_h^+ \dd s \nno \\
&&   - \int_{\extbdr_D} 
   \left(G_{\Gamma}^{\top}({\bf u}_h^+) \nabla_h{\bf v}^+_h\right) :
\left({\bf u}_h^+-{\bf  u}_\Gamma({\bf u}_h^+) \right) \otimes {\bf n} \dd s, \nno
\end{eqnarray}
where $\uuu{\delta}_{\Gamma}({\bf u}_h) = \Cip \tfrac{\polyo^2}{\hh}
G_\Gamma({\bf u}_h^+)\left({\bf u}_h -{\bf u}_\Gamma({\bf u}_h)\right) \otimes
{\bf n}$. Here, the viscous boundary flux ${\mathcal F}_{\Gamma}^{v}$ and the
corresponding homogeneity tensor $G_{\Gamma}$  are defined by
\begin{equation*}
{\mathcal F}_{\Gamma}^{v}({\bf u}_h, \nabla {\bf u}_h)
= {\mathcal F}^{v}({\bf u}_{\Gamma}({\bf u}_h), \nabla {\bf u}_h)
= G_{\Gamma}({\bf u}_h)\nabla {\bf u}_h = G({\bf u}_{\Gamma}({\bf u}_h)) \nabla {\bf u}_h.
\end{equation*}
The convective boundary
flux ${\mathcal H}_\Gamma$ is defined by
\begin{equation*}
  {\mathcal H}_\Gamma({\bf u}_h^+,{\bf u}_\Gamma({\bf u}_h^+),{\bf n})
={\bf n}\cdot{\mathcal F}^c({\bf u}_\Gamma({\bf u}_h^+)).
\end{equation*}
Finally, the boundary function ${\bf u}_\Gamma({\bf u})$ is given according to
the type of boundary condition imposed; in the current setting ${\bf
u}_\Gamma({\bf u}) = {\bf g}_D$ on $\extbdr_D$ and ${\bf u}_\Gamma({\bf u}) =
{\bf u}$ on $\extbdr_N$. For further details regarding the imposition
of more general boundary conditions, we refer to~\cite{HH08a}, for
example.

\section{Computational Framework for DGFEMs}
\label{sec:comp_framework_for_dgfems}

In this section we present the general computational framework for the
automatic generation of DGFEM (semi-) linear forms in a concise
and easy to use manner. We begin by outlining the treatment of 
both convective and viscous
components arising in conservation laws. In this setting, the DGFEM
formulations can be constructed in a consistent manner. We exploit this by
designing utility functions to automatically generate these symbolic DGFEM
formulations. We then proceed by proposing a hierarchical framework for computing
DGFEM formulations of PDE operators. This hierarchy takes advantage of the
DGFEM utility functions, providing a modular framework for a suite of DGFEM
formulations for various operators arising in PDE problems of engineering
interest.

%
\subsection{Automatic Treatment of Convective Terms}
\label{sec:convective_auto_dg}
In this section we discuss the automatic symbolic representation of the DGFEM
discretisation of the term involving the convective flux function
$\mathcal{F}^c\left(\cdot\right)$. To this end, we recall that the convective
numerical flux function $\mathcal{H}(\cdot, \cdot, \cdot)$, cf.~\eqref{eq:dg-ns}, 
must be defined on the element interfaces. Defining a callable
function \ilc{F_c} which specifies $\mathcal{F}^c\left(\cdot\right)$ we
construct the abstraction of the numerical flux function as shown in
Listing~\ref{code:abstract_convective_flux}. The methods \ilc{interior} and
\ilc{exterior} will be called to construct the flux function on interior and
exterior faces, respectively. The method \ilc{setup} is provided for the 
inheriting class to initialise any members prior to calls to \ilc{interior}
and \ilc{exterior}. This design allows for the flux function to be different on 
the two types of faces present in the mesh.

\begin{table}[t!]
    \begin{lstlisting}[language=Python, caption={The abstraction of the
        numerical flux function $\mathcal{H}(\cdot, \cdot, \cdot)$.},captionpos=b , label=code:abstract_convective_flux]
class ConvectiveFlux:

    def __init__(self):
        pass

    def setup(self):
        pass

    def interior(self, F_c, u_p, u_m, n):
        pass

    def exterior(self, F_c, u_p, u_m, n):
        pass
    \end{lstlisting}
\end{table}

For simplicity, here we consider two implementations of the
\ilc{ConvectiveFlux} class, based on employing the local Lax-Friedrichs and HLLE fluxes,
though we stress that other fluxes, such as the Vijayasundaram or Roe
flux, for example, may also be employed. Thereby, on the boundary
$\partial\kappa$ of an element $\kappa\in\mesh$, we define the local-Lax
Friedrichs flux $\mathcal{H}_\text{LF}$ and HLLE flux $\mathcal{H}_\text{HLLE}$ by
%
%
\begin{subequations}
\begin{align}
\mathcal{H}_\text{LF}\left.\left(\mathbf{u}^+_h, \mathbf{u}^-_h, \mathbf{n}_\kappa\right)\right|_{\partial\kappa}
&= \tfrac{1}{2}\left( \mathcal{F}^c\left(\mathbf{u}^+_h\right) \cdot \mathbf{n}_\kappa + \mathcal{F}^c\left(\mathbf{u}^-_h\right) \cdot \mathbf{n}_\kappa + \alpha \left(\mathbf{u}^+_h - \mathbf{u}^-_h\right) \right), \label{eq:local_lf_flux} \\
\mathcal{H}_\text{HLLE}\left.\left(\mathbf{u}^+_h, \mathbf{u}^-_h, \mathbf{n}_\kappa\right)\right|_{\partial\kappa}
&= \tfrac{1}{\lambda^+ - \lambda^-}
\lambda^+ \mathcal{F}^c(\mathbf{u}^+_h) \cdot \mathbf{n}_\kappa
-\lambda^- \mathcal{F}^c(\mathbf{u}^-_h) \cdot \mathbf{n}_\kappa
-\lambda^+ \lambda^- (\mathbf{u}^+_h - \mathbf{u}^-_h),
\end{align}
\end{subequations}
respectively. Here, $\alpha$ is the local dissipation parameter, which is selected to be the
maximum of the (absolute value) of the eigenvalues of the Jacobi matrix
\begin{equation*}
B(\mathbf{u},\mathbf{n}_\kappa)
= \sum_{i=1}^d \frac{\partial \mathbf{f}^c_i}{\partial\mathbf{u}}
\mathbf{n}_{\kappa,i},
\end{equation*}
evaluated on the element face. More precisely, we have that
$$\left.\alpha\right|_{\partial\kappa} = \max_{\mathbf{w}=\mathbf{u}^+_h,\mathbf{u}^-_h} \left \lbrace | \lambda_{\max}
\left( B \left( \mathbf{w}, \mathbf{n}_\kappa \right) \right )  | \right
\rbrace$$, where $\lambda_{\max}$ denotes the largest eigenvalue (in modulus) of 
$B \left( \cdot, \mathbf{n}_\kappa \right)$. Additionally,
\begin{eqnarray*}
\lambda^+ &=& \max \left( \max_{\mathbf{w}=\mathbf{u}^+_h,\mathbf{u}^-_h} \left \lbrace \lambda_{\max}
\left( B \left( \mathbf{w}, \mathbf{n}_\kappa \right) \right ) \right \rbrace, 0 \right), 
\nonumber \\
\lambda^- &=& \min \left( \min_{\mathbf{w}=\mathbf{u}^+_h,\mathbf{u}^-_h} \left \lbrace \lambda_{\min}
\left( B \left( \mathbf{w}, \mathbf{n}_\kappa \right) \right ) \right \rbrace, 0 \right), 
\end{eqnarray*}
where $\lambda_{\min}$ denotes the smallest eigenvalue (in modulus) of 
$B \left( \cdot, \mathbf{n}_\kappa \right)$.

Given \ilc{F_c}, which specifies
$\mathcal{F}^c\left(\cdot\right)$, the numerical flux functions
$\mathcal{H}_\text{LF}(\cdot,\cdot,\cdot)$ and
$\mathcal{H}_\text{HLLE}(\cdot,\cdot,\cdot)$ can be automatically generated in
order to yield the discretisation of the convective term present in the
underlying PDE problem; we note that the constructors of both classes
\ilc{LocalLaxFriedrichs} and \ilc{HLLE} require the 
symbolic representation of the eigenvalues of $B$.
%
%
As an example, we consider the application of the DGFEM employing the local Lax-Friedrichs flux
for the numerical approximation of the linear advection equation
\begin{equation}
	\label{eq:linear_advection_ufl}
	\nabla \cdot \left( \mathbf{b} u \right) = f,
\end{equation}
where $\mathbf{b} =(b_1,b_2,\ldots,b_d)^\top$, 
$\mathcal{F}^c(u) = \mathbf{b} u$, and $f$ is some given forcing function. 
Here, the dissipation parameter can be shown to be
$\left.\alpha\right|_{\partial\kappa} = \left| \mathbf{b} \cdot
\mathbf{n}_\kappa \right|$. The UFL code required to generate the numerical 
flux function, together with the corresponding element boundary term arising 
in the DGFEM scheme for this problem is given in
Listing~\ref{code:linear_advection_autogen}. We note that 
the use of the class \ilc{HLLE} defining the HLLE numerical flux follows in an 
analogous manner. 

\begin{table}[t!]
	\begin{lstlisting}[language=Python, caption={Example of the automatic calculation and symbolic representation of the local Lax-Friedrichs flux $\mathcal{H}_\text{LF}(u^+, u^-,\mathbf{n}_\kappa)$ for the linear advection equation shown in~\eqref{eq:linear_advection_ufl}.},captionpos=b , label=code:linear_advection_autogen]
def F_c(u): return b*u
H = LocalLaxFriedrichs(lambda u, n: dot(b, n))
H.setup(F_c, u('+'), u('-'), n('+'))
conv_interior = H.interior(F_c, u('+'), u('-'), n('+'))*(v('+') - v('-'))*dS
	\end{lstlisting}
\end{table}

In many cases an analytical expression for the eigenvalues of the Jacobi matrix
$B(\mathbf{u},\mathbf{n}_\kappa)$, $\kappa\in\mesh$,
may be difficult to evaluate; thereby, packages such as SymPy~\cite{sympy} may be used to
compute them symbolically. To this end, $B(\mathbf{u},
\mathbf{n}_\kappa)$ can be assembled using symbolic differentiation; the
eigenvalues of this matrix may then be computed symbolically by exploiting the Berkowitz
algorithm~\cite{BERKOWITZ1984147}. In Listing~\ref{code:eigenvaluesalpha} we show 
an example of this method applied
to the convective component of the incompressible Navier-Stokes equations,
i.e., $\mathcal{F}^c(\cvec{u}) = \cvec{u} \otimes
\cvec{u}$.

\begin{table}[t!]
  \begin{lstlisting}[language=Python, caption={Automatic symbolic algebra computation of 
    flux Jacobian eigenvalues
  required by the local Lax-Friedrichs and HLLE fluxes applied to the incompressible Navier-
  Stokes convective flux component $\mathcal{F}^c(\cvec{u}) = \cvec{u}
  \otimes
  \cvec{u}$.},captionpos=b , label=code:eigenvaluesalpha]
from sympy import *

dim = 2
u = Matrix([Symbol("u%d" % d, real=True) for d in range(dim)])
n = Matrix([Symbol("n%d" % d, real=True) for d in range(dim)])

F_c = u*u.T

B = zeros(dim, dim)
for d in range(dim):
  B += F_c[:, d].jacobian(u)*n[d]

print(B.eigenvals().keys())
  \end{lstlisting}
\end{table}

\subsection{Automatic Treatment of Viscous Terms}
\label{subsec:automatic_treatment_of_viscous_terms}
In this section we develop utility functions which
automatically generate the DGFEM discretisation of
second--order PDE operators. To this end, we recall from
Section~\ref{sec:dg-fems} that the viscous component of the
underlying PDE is given by
\begin{equation}
	- \nabla \cdot \mathcal{F}^v\left(\mathbf{u}, \nabla \mathbf{u}\right) = \mathbf{0};
    \label{eq:visc_term_only}
\end{equation}
here, we have set the right-hand in \eqref{eq:visc_term_only} to zero, for simplicity,
in order to concentrate on the discretisation of the viscous term.
The semilinear DGFEM formulation of~\eqref{eq:visc_term_only}, shown in
equation~\eqref{eq:dg-ns}, is encapsulated by implementations of the
abstract class
$$
\mbox{\ilc{DGFemViscousTerm(F_v, u, v, gamma, G, n)}} ,
$$
where 
\begin{equation}
	\text{\ilc{F_v(u, grad\_u)}} = \mathcal{F}^v(\mathbf{u}, \nabla \mathbf{u})
\end{equation}
is a callable function, which defines the viscous flux function,
\ilc{u} and \ilc{v} denote the trial and test functions, respectively,
$\text{\ilc{gamma}} \equiv \gamma = \Cip \nicefrac{\polyo^2}{\hh}$
is the DGFEM penalisation coefficient, $\text{\ilc{G}}\equiv G(\cvec{u})$ is the homogeneity tensor, and $\text{\ilc{n}}\equiv \cvec{n}$
is the face normal.



\begin{table}[t!]
    \begin{lstlisting}[language=Python, caption={The class \ilc{DGFemViscousTerm} provides the 
    abstract interface for interior and exterior residual formulation using
    symbolic algebra. The class \ilc{DGFemSIPG} uses UFL to automatically
    formulate the interior and exterior terms of the SIPG 
    formulation of the viscous component of~\eqref{eq:dg-ns}.}, 
    captionpos=b, label=code:dg_visc_term_class] 
class DGFemSIPG(DGFemViscousTerm):

  def interior_residual(self, dInt):
    G = self.G
    F_v, u, v, grad_v = self.F_v, self.U, self.V, self.grad_v_vec
    sig, n = self.sig, self.n

    residual = \
      - inner(tensor_jump(u, n), avg(hyper_tensor_T_product(G, grad_v)))*dInt \
      - inner(ufl_adhere_transpose(avg(F_v(U))), tensor_jump(v, n))*dInt \
      + inner(sig('+')*hyper_tensor_product(g_avg(G), tensor_jump(u, n)), tensor_jump(v, n))*dInt
    return residual

  def exterior_residual(self, u_gamma, dExt):
    G = self._make_boundary_G(self.G, u_gamma)
    F_v, u, v, grad_u, grad_v = self.F_v, self.U, self.V, grad(self.U), self.grad_v_vec
    n = self.n

    residual = \
      - inner(dg_outer(u - u_gamma, n), hyper_tensor_T_product(G, grad_v)) * dExt \
      - inner(hyper_tensor_product(G, grad_u), dg_outer(v, n)) * dExt \
      + inner(self.sig*hyper_tensor_product(G, dg_outer(u - u_gamma, n)), dg_outer(v, n)) * dExt
    return residual
    \end{lstlisting}
\end{table}

Recalling the definition of the homogeneity tensor~\eqref{eq:homo_tensor} and
the homogeneity tensor product~\eqref{eq:homo_tens_prod}, $G(\cvec{u})$ is
automatically computed using the function \ilc{homogeneity_tensor(F_v, u)}
and the function
\ilc{hyper_tensor_product(G, tau)} computes the homogeneity tensor product.
These functions may then be
employed for the generation of the corresponding DGFEM formulation. On the basis
of these two functions, we introduce the abstract class
\ilc{DGFemViscousTerm}; this offers the following three methods for handling
each of the boundary components present in the DGFEM scheme:
\begin{enumerate}
    \item  \ilc{DGFemViscousTerm.interior_residual(dS)} automatically generates
        terms associated with the interior boundaries $\intbdr$ present
        in~\eqref{eq:dg-ns};
    \item \ilc{DGFemViscousTerm.exterior_residual(u_gamma, ds_i)}
        generates the terms associated with exterior boundary component
        \ilc{ds_i} present in~\eqref{eq:dg-ns-exterior} with boundary condition
        $\mathbf{u}_\Gamma(\mathbf{u}) = \text{\ilc{u_gamma}}$;
    \item Finally, \ilc{DGFemViscousTerm.neumann_residual(gN, ds_i)} generates
        any terms arising from Neumann boundary conditions with flux
        specification $\mathcal{F}^v(\mathbf{u}, \nabla \mathbf{u}) \cdot
        \mathbf{n} =
        \cvec{g}_N$.
\end{enumerate}
We demonstrate an example of the implementation of the
\ilc{DGFemViscousTerm} with the SIPG method in
Listing~\ref{code:dg_visc_term_class}; for the sake of brevity, we have removed
input and error checking. Further implementations such as the non--symmetric
interior penalty (NIPG) and Baumann--Oden schemes are available in the classes
\ilc{DGFemNIPG} and \ilc{DGFemBO}, respectively. The DGFEM formulation 
proposed by Bassi \& Rebay \cite{BRMPS97} is challenging in the symbolic 
framework due to the requirement of solving
elementwise problems for the local lifting operator. Such operations are not
easily formulated in the UFL for use with DOLFIN, for example; the 
implementation of DGFEMs defined based on employing local lifting operators will
be considered as part of our programme of future research. For a more detailed
discussion of the formulation of lifting operators in DOLFIN, we refer to~\cite[Chapter
5]{oelgaard_phd}.

An example of the UFL code required to generate the DGFEM semilinear residual
discretisation of the quasi-linear second-order PDE problem
\begin{subequations}
\label{eq:quick_poisson_for_dg_visc}
\begin{align}
    -\nabla \cdot \left((u + 1) \nabla u \right) & = f \; \text{in} \; \Omega, \\
    u & = g_D \; \text{on} \; \Gamma_D,  \\
    \nabla u \cdot \cvec{n} & = g_N \; \text{on} \; \Gamma_N 
\end{align}
\end{subequations}
is shown in Listing~\ref{code:full_dgfemvisc_poi_ex}.

\begin{table}[t!]
	\begin{lstlisting}[language=Python,caption={Example UFL code for the DGFEM discretisation of the quasi-linear PDE~\eqref{eq:quick_poisson_for_dg_visc} using the
    \ilc{DGFemSIPG} utility class.}, captionpos=b,label=code:full_dgfemvisc_poi_ex]
def F_v(u, grad_u): return (u + 1)*grad(u)
G = homogeneity_tensor(F_v, u)
vt = DGFemSIPG(F_v, u, v, sig, G, n)
residual = dot((u + 1)*grad(u), grad(v))*dx \
         + vt.interior_residual(dS) \
         + vt.exterior_residual(g_D, ds_D) \
         + vt.neumann_residual(g_N, ds_N) \
         - f*v*dx
	\end{lstlisting}
\end{table}

\subsection{Automatic Generation of DGFEM Formulations}
Even with the utility functions outlined in Sections
\ref{sec:convective_auto_dg} and
\ref{subsec:automatic_treatment_of_viscous_terms}, the specification of the
DGFEM scheme can be very verbose; this is particularly the case when
discretising systems comprising many PDE variables with multiple boundary
conditions. In this section we propose a hierarchical scheme for the management
of DGFEM formulations of PDE operators of increasing complexity.

\subsubsection{Boundary Conditions}

Firstly, we outline a simple framework for managing the boundary conditions
enforced within a given PDE problem. We define the \ilc{DGBC} (discontinuous
Galerkin boundary condition) abstract class from which the classes
\ilc{DGDirichletBC} and \ilc{DGNeumannBC} inherit. These implementations simply
serve to store the boundary condition and the boundary component over which the
condition should be enforced. For example, applying a Dirichlet boundary
condition as required by the quasi-linear PDE problem stated
in~\eqref{eq:quick_poisson_for_dg_visc} simply requires instantiation of
\ilc{DGDirichletBC(ds\_D, g\_D)}. Similarly, for the imposition of
the Neumann boundary condition we construct
\ilc{DGNeumannBC(ds\_N, g\_N)}; note that Robin conditions can be implemented
in an analogous manner. By generating a list of boundary conditions in
this manner, we may easily formulate their imposition within the DGFEM scheme.
As an example of a series of Dirichlet boundary conditions being automatically
generated using \ilc{DGFemViscousTerm}, we refer to
Listing~\ref{code:auto_dirichlet_dg}.

\begin{table}[t!]
	\begin{lstlisting}[language=Python,caption={Example of the automatic generation of the exterior residual terms of a given \ilc{DGFemViscousTerm} for a list of Dirichlet boundary conditions.}, captionpos=b,
	label=code:auto_dirichlet_dg]
bcs = [DGDirichletBC(bc_1, ds_1), DGDirichletBC(bc_2, ds_2), ...]
vt = DGFemViscousTerm(F_v, u, v, delta)
ext = sum(vt.exterior_residual(bc.get_function(), bc.get_boundary()) for bc in bcs)
	\end{lstlisting}
\end{table}  

\subsubsection{Abstract DGFEM Formulation}

\begin{table}[t!]
    \begin{lstlisting}[language=Python, caption={The \ilc{HyperbolicOperator} class.},captionpos=b,
    label=code:hyperbolic_op_class]
class HyperbolicOperator(DGFemFormulation):

    def __init__(self, mesh, V, bcs,
                 F_c=lambda u: u,
                 H=LocalLaxFriedrichs(lambda u, n: inner(u, n))):
        DGFemFormulation.__init__(self, mesh, V, bcs)
        self.F_c = F_c
        self.H = H

    def generate_fem_formulation(self, u, v, dx=None, dS=None):
        if dx is None:
            dx = Measure('dx', domain=self.mesh)
        if dS is None:
            dS = Measure('dS', domain=self.mesh)
        n = FacetNormal(self.mesh)

        residual = -inner(self.F_c(u), grad(v))*dx

        self.H.setup(self.F_c, u('+'), u('-'), n('+'))
        residual += inner(self.H.interior(self.F_c, u('+'), u('-'), n('+')), (v('+') - v('-')))*dS

        for bc in self.dirichlet_bcs:
            gD = bc.get_function()
            dSD = bc.get_boundary()

            self.H.setup(self.F_c, u, gD, n)
            residual += inner(self.H.exterior(self.F_c, u, gD, n), v)*dSD

        for bc in self.neumann_bcs:
            dSN = bc.get_boundary()

            residual += inner(dot(self.F_c(u), n), v)*dSN

        return residual
    \end{lstlisting}
\end{table} 

We encapsulate the abstraction of a DGFEM scheme in the class
\ilc{DGFemFormulation} which prescribes one abstract method
\ilc{generate\_fem\_formulation.} The \ilc{DGFemFormulation} constructor
requires the mesh $\omegatri$, the function space $\mathbf{V}^m_\polyo$ and the
list of boundary conditions. This class will serve as the base class for the
DGFEM formulation of all derived PDE operators. In this work we describe two
direct children of \ilc{DGFemFormulation}: \ilc{HyperbolicOperator} and
\ilc{EllipticOperator}. We stress that the flexibility of this design permits
DGFEM formulations of other PDE operators, cf. Sections~\ref{sec:maxwell},
\ref{sec:hyperelasticity}, and \ref{sec:hyperelasticity2} below.

The class \ilc{HyperbolicOperator} inherits \ilc{DGFemFormulation} and its purpose
is to generate DGFEM formulations of the PDE operator $\nabla \cdot
\mathcal{F}^c(\cdot)$. The class implementation is shown in
Listing~\ref{code:hyperbolic_op_class}; here, \ilc{generate\_fem\_formulation}
is overridden to construct the DGFEM formulation of the provided convective flux
operator on the interior and exterior faces, as well as on the elements in the mesh. The numerical
flux function provided in the constructor \ilc{H}, which implements
\ilc{ConvectiveFlux}, is used for the DGFEM formulation.
%
%
To highlight to modularity of this design, consider an extension of the
\ilc{HyperbolicOperator} for the time-dependent Burgers' equation
\begin{equation}
\frac{\partial u}{\partial t} + \frac{\partial}{\partial x} (\nicefrac{1}{2} \, u^2) = 0.
\label{eq:burgers}
\end{equation}
Setting $t=y$, we may recast \eqref{eq:burgers} in the following equivalent form
\begin{equation}
    \nabla \cdot \mathcal{F}^c(u) \equiv
    \nabla \cdot
    \begin{pmatrix}
    \frac{1}{2} u^2 \\
    u
    \end{pmatrix}
    = 0;
\end{equation}{}
the implementation of the class \ilc{SpacetimeBurgersOperator} (using
\ilc{LocalLaxFriedrichs} by default) is depicted in
Listing~\ref{code:burgers_op_class}. Here, we note that the derived class 
must simply specify the form of the convective flux $\mathcal{F}^c(u)$ and the
flux function $\mathcal{H}(\cdot,\cdot,\cdot)$; indeed, once these are defined,
the automatic generation of the DGFEM formulation is then handled by the parent
\ilc{DGFemFormulation}.

\begin{table}[t!]
    \begin{lstlisting}[language=Python, caption={The \ilc{SpacetimeBurgersOperator} class.}, captionpos=b,
    label=code:burgers_op_class]
class SpacetimeBurgersOperator(HyperbolicOperator):

    def __init__(self, mesh, V, bcs, flux=None):

        def F_c(u):
            return as_vector((u**2/2, u))

        if flux is None:
            flux = LocalLaxFriedrichs(lambda u, n: u*n[0] + n[1])

        HyperbolicOperator.__init__(self, mesh, V, bcs, F_c, flux)
    \end{lstlisting}
\end{table} 

Secondly, we introduce the class \ilc{EllipticOperator}, which inherits
\ilc{DGFemFor\-mulation}, and requires the specification of
$\mathcal{F}^v(\cdot, \nabla \cdot)$ at instantiation. The overridden method
\ilc{generate\_fem\_formulation} is then written to automatically generate the
interior and exterior boundary, as well as the element, integration terms, implementing
all of the concepts of the utility functions for elliptic operators in
\ilc{DGFemViscousTerm}. The outline of the \ilc{EllipticOperator} class,
which generates the SIPG formulation by default, is given in
Listing~\ref{code:elliptic_op_class}. To highlight the modularity of
\ilc{EllipticOperator} in Listing~\ref{code:poisson_op_class} we show the
implementation of the class \ilc{PoissonOperator} which specifies
$\mathcal{F}^v(u,\nabla u) = {\mathcal K} \nabla u$, where ${\mathcal K}$ is the
diffusion coefficient. An example of the automatic generation of the DGFEM
formulation of the quasi-linear elliptic PDE problem
\eqref{eq:quick_poisson_for_dg_visc} is given in
Listing~\ref{code:poisson_auto_gen_dg}.

\begin{table}[t!]
    \begin{lstlisting}[language=Python,caption={The \ilc{EllipticOperator} class.},captionpos=b,
    label=code:elliptic_op_class]
class EllipticOperator(DGFemFormulation):

    def __init__(self, mesh, fspace, bcs, F_v, C_IP=10.0):
        DGFemFormulation.__init__(self, mesh, fspace, bcs)
        self.F_v = F_v
        self.C_IP = C_IP

    def generate_fem_formulation(self, u, v, dx=None, dS=None, vt=None):
        if dx is None:
            dx = Measure('dx', domain=self.mesh)
        if dS is None:
            dS = Measure('dS', domain=self.mesh)

        h = CellVolume(self.mesh)/FacetArea(self.mesh)
        n = FacetNormal(self.mesh)
        sigma = self.C_IP*Constant(max(self.fspace.ufl_element().degree()**2, 1))/h
        G = homogeneity_tensor(self.F_v, u)

        if vt is None:
            vt = DGFemSIPG(self.F_v, u, v, sigma, G, n)

        if inspect.isclass(vt):
            vt = vt(self.F_v, u, v, sigma, G, n)

        assert(isinstance(vt, DGFemViscousTerm))

        residual = inner(self.F_v(u, grad(u)), grad(v))*dx
        residual += vt.interior_residual(dS)

        for dbc in self.dirichlet_bcs:
            residual += vt.exterior_residual(dbc.get_function(), dbc.get_boundary())

        for dbc in self.neumann_bcs:
            residual += vt.neumann_residual(dbc.get_function(), dbc.get_boundary())

        return residual
    \end{lstlisting}
\end{table}

\subsubsection{Hierarchy of PDE Operators}

A natural extension of this framework is a hierarchy of automatically 
generated DGFEM operators. As shown in Listing~\ref{code:burgers_op_class},
the \ilc{SpacetimeBurgersOperator} inherits from \ilc{HyperbolicOperator}.
Consider now, for example, the compressible Navier Stokes equations.
Here the inheritance chain may begin with the \ilc{HyperbolicOperator}, from
which a \ilc{CompressibleEulerOperator} would inherit, and in turn a 
\ilc{CompressibleNavierStokesOperator} would inherit \ilc{Compress\-ibleEulerOperator}
and additionally \ilc{EllipticOperator} for the viscosity terms. Further implementations
of each member of this class hierarchy need not only be undertaken by
inheritance. Operators deriving from models of large physical systems may more
appropriately aggregate sub-operators as necessary. The derived classes must
simply manage the function spaces and boundary conditions amongst the
aggregated \ilc{DGFemFormulation} members. This framework significantly
reduces the code required for subsequent development of DGFEM formulations of
PDE operators of increasing complexity. By ensuring that each layer of the
hierarchy is correctly verified and fully tested means that DGFEM formulations
may be debugged in a very straightforward manner.

\begin{table}[t!]
    \begin{lstlisting}[language=Python,caption={The \ilc{PoissonOperator} class need only inherit the \ilc{EllipticOperator}
    and define its own viscous flux, \ilc{F\_v}.},captionpos=b,
    label=code:poisson_op_class]
class PoissonOperator(EllipticOperator):

    def __init__(self, mesh, fspace, bcs, kappa=1):
        def F_v(u, grad_u):
            return kappa*grad_u

        EllipticOperator.__init__(self, mesh, fspace, bcs, F_v)
    \end{lstlisting}
\end{table}

\begin{table}[t!]
	\begin{lstlisting}[language=Python, caption={Example of implementing the \ilc{PoissonOperator} utility class.},captionpos=b,
	label=code:poisson_auto_gen_dg]
po = PoissonOperator(mesh, V, DGDirichletBC(ds, gD), kappa=u+1)
residual = po.generate_fem_formulation(u, v) - f*v*dx
	\end{lstlisting}
\end{table}

\section{Examples} \label{sec:examples}

In this section we present a series of examples of increasing complexity
to highlight the ease in which each PDE problem may be discretised using
a DGFEM formulation, based on employing the class hierarchy proposed
and implemented within this article. We stress that the verbosity and
complexity of specifying each DGFEM discretisation is vastly reduced 
within this modular framework. To this end, we consider the
discretisation of a simple scalar advection-diffusion equation, the 
compressible Euler equations, the compressible Navier-Stokes equations posed
in both conserved and entropy variables, the indefinite Maxwell problem, 
and a hyperelasticity problem.
In each case, for simplicity of presentation, we employ the SIPG 
DGFEM discretisation of the second-order PDE operator and a local Lax-Friedrichs flux 
for the numerical approximation of the convective terms.
For brevity, in some of the examples given below only code snippets will
be shown, however, full versions of the corresponding python scripts are
available from \url{https://bitbucket.org/nate-sime/dolfin_dg}.

\subsection{Example 1: Advection-diffusion problem}
\label{sec:diffconvreac}

In this first example, we highlight the use of the classes
\ilc{HyperbolicOperator} and \ilc{EllipticOperator} given in
Listings~\ref{code:hyperbolic_op_class} \&~\ref{code:elliptic_op_class},
respectively. To this end, given $\domain = (0, 1)^2$, with boundary $\extbdr$, 
consider the problem: find
$u$ such that
\begin{align}
\label{eq:diffconvreac_apriori}
-\nabla \cdot \left( {\mathcal K} \nabla u\right) + \nabla \cdot \left( \cvec{b} u^2 \right) &= f  &\text{in} \; \domain, \\
u &= g_D &\text{on} \; \extbdr, \label{eq:diffconvreac_apriori_bc}
\end{align}
where ${\mathcal K} > 0$ denotes the diffusion coefficient 
and $\cvec{b} = (b_1,b_2)^\top$. Thereby, setting $m=1$ and $d=2$,
the PDE problem~\eqref{eq:diffconvreac_apriori}, \eqref{eq:diffconvreac_apriori_bc} 
can be written in the general
form~\eqref{eq:model-prob}--\eqref{eq:bcs2}, where $\mathcal{F}^c(\mathbf{u}) =
\mathbf{b} u^2$,
$\mathcal{F}^v(\mathbf{u}, \nabla{\mathbf{u}}) = {\mathcal K} \nabla u$,
and $\extbdr[N] = \emptyset$. Moreover, it can easily be shown that 
the dissipation parameter arising in the local Lax-Friedrichs
flux,  is given by 
$\left.\alpha\right|_{\partial\kappa} =
\max_{w = u^+, u^-} 2\left| w \cvec{b} \cdot \mathbf{n}_\kappa \right|$,
while the homogeneity tensor $G$, cf.~\eqref{eq:homo_tensor},
has entries $G_{11} = G_{22} = {\mathcal K}$, $G_{12} = G_{21} = 0$.

\begin{figure}[t!]
\centering
\begin{subfigure}{.5\textwidth}
  \centering
  \includegraphics[width=1.\linewidth]{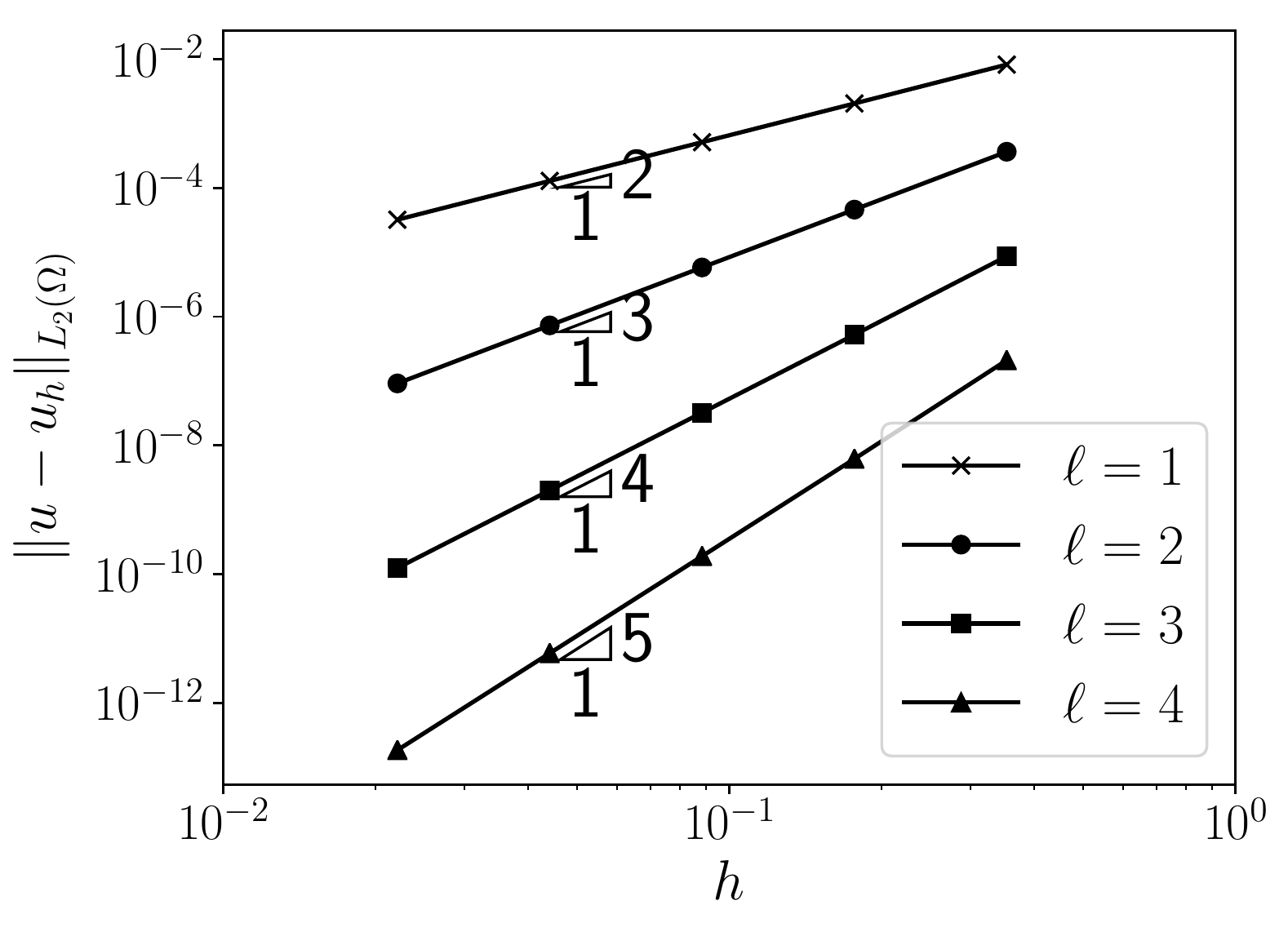} \\
  (a)
\end{subfigure}%
\begin{subfigure}{.5\textwidth}
  \centering
  \includegraphics[width=1.\linewidth]{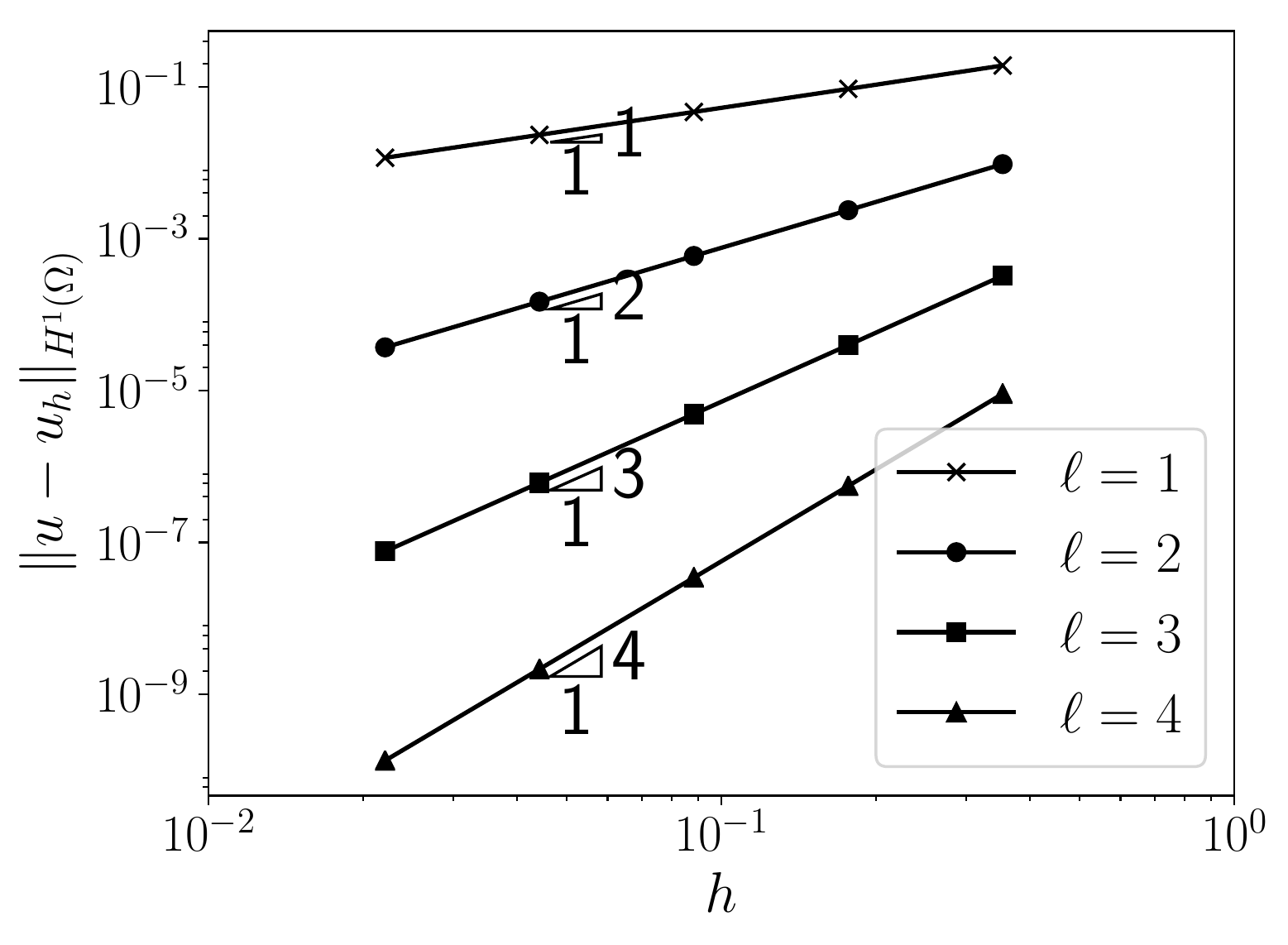} \\
  (b)
\end{subfigure}
\caption{Example 1: Convergence of the DGFEM with $h$--refinement: 
(a) $\|u-u_h\|_{L_2(\Omega)}$; (b) $\|u-u_h\|_{H^1(\Omega)}$.}
\label{fig:ad_conv_rates}
\end{figure}

The specification of the numerical fluxes $\mathcal{F}^c$ and $\mathcal{F}^v$,
when ${\mathcal K} = 1+u$, along with $\alpha(u, \cvec{n})$ used with the
local Lax-Friedrichs flux function within UFL syntax is provided in the
code listing presented in Listing~\ref{code:ufl_ad_fcfv}. Exploiting the
software concepts outlined in Section~\ref{sec:comp_framework_for_dgfems} for
the automatic computation of DGFEM formulations, the UFL code to construct and
solve the resulting DGFEM approximation of ~\eqref{eq:diffconvreac_apriori},
\eqref{eq:diffconvreac_apriori_bc} is presented in
Listing~\ref{code:ufl_ad_full}, where we have specified that $\cvec{b} =
(1, 1)^\top$, $f=-4\mbox{{\rm e}}^{2(x - y)} - 2\mbox{{\rm e}}^{x - y}$, and
$g_D = e^{x - y}$; thereby, the analytical solution
to~\eqref{eq:diffconvreac_apriori}, \eqref{eq:diffconvreac_apriori_bc} is given
by $u = e^{x - y}$. We note that on the final line of the code listing given in
Listing~\ref{code:ufl_ad_full}, the call to DOLFIN's
\ilc{solve} function automatically computes the G\^{a}teaux derivative
of the DGFEM semilinear form ${\mathcal N}(\cdot,\cdot)$, cf. ~\eqref{eq:dg-ns},
employed for the numerical approximation of~\eqref{eq:diffconvreac_apriori},
\eqref{eq:diffconvreac_apriori_bc} and utilises a Newton solver to evaluate the
DGFEM solution; for further details, we refer to~\cite{Simula.simula.1092}.
Finally, in Figure~\ref{fig:ad_conv_rates} we show the asymptotic behaviour of
the underlying DGFEM on a sequence of uniformly refined triangular meshes with
polynomial orders $\ell=1,2,3,4$; as we expect, we observe that the
$L_2(\Omega)$ and $H^1(\Omega)$ norms of the error tend to zero at the
respective rates of ${\mathcal O}(h^{\ell+1})$ and ${\mathcal O}(h^{\ell})$ as
the mesh size $h$ tends to zero. The complete code employed for this numerical
example is provided in the file \url{advection_diffusion.py}.


\begin{table}[t!]
\begin{lstlisting}[language=Python, caption={Example 1: UFL representation of $\mathcal{F}^c$ and 
$\mathcal{F}^v$ of~\eqref{eq:diffconvreac_apriori}.},captionpos=b, label=code:ufl_ad_fcfv]
    def F_c(u):
        return b*u**2

    H = LocalLaxFriedrichs(lambda u, n: 2*u*dot(b, n))

    def F_v(u, grad_u):
        return (u + 1)*grad_u
\end{lstlisting}
\end{table}
\begin{table}[t!]
\begin{lstlisting}[language=Python,caption={Example 1: Automatic DGFEM formulation for the numerical approximation of~\eqref{eq:diffconvreac_apriori}, \eqref{eq:diffconvreac_apriori_bc}, using the definitions of the
    fluxes in Listing~\ref{code:ufl_ad_fcfv}.},captionpos={}b,
label=code:ufl_ad_full]
mesh = UnitSquareMesh(32, 32)

V = FunctionSpace(mesh, 'DG', 1)
u, v = Function(V), TestFunction(V)

gD = Expression('exp(x[0] - x[1])', element=V.ufl_element())
f = Expression('-4*exp(2*(x[0] - x[1])) - 2*exp(x[0] - x[1])',
               element=V.ufl_element())
b = Constant((1, 1))

ho = HyperbolicOperator(mesh, V, DGDirichletBC(ds, gD), F_c, H)
eo = EllipticOperator(mesh, V, DGDirichletBC(ds, gD), F_v)

residual = ho.generate_fem_formulation(u, v) \
           + eo.generate_fem_formulation(u, v) \
           - f*v*dx

solve(residual == 0, u)
\end{lstlisting}
\end{table}


\subsection{Example 2: Compressible Euler Equations}

In this second example, we consider the DGFEM discretisation of
the compressible Euler equations: find $\cvec{u} = (\rho,
\rho\cvec{u}_{\text{v}}, \rho E)^\top : \domain \rightarrow \mathbb{R}^+ \times
\mathbb{R}^{d} \times \mathbb{R}^+$ such that
\begin{alignat}{2}
  \nabla \cdot \mathcal{F}^c \left( \cvec{u} \right) &= \cvec{0} &&\quad\text{in } \domain, \label{eq:comp_euler_bvp} 
\end{alignat}
where the convective flux is given by
\begin{equation}
  \mathcal{F}^c(\cvec{u}) = 
  \begin{pmatrix}
    \rho \cvec{u}_{\text{v}} \\
    \rho \cvec{u}_{\text{v}} \otimes \cvec{u}_{\text{v}} + p \eyedentity \\
    \left( \rho E + p\right) \cvec{u}_{\text{v}} \\
  \end{pmatrix}.
  \label{eq:euler_flux}
\end{equation}
Here, $\rho$, ${\bf u}_{\text{v}}$, $p$, and $E$ denote the density,
velocity vector, pressure, and specific total energy, respectively.
The equation of state of an ideal gas is given by
$$
p = (\gamma -1) \rho \left( E - \frac{1}{2} |{\bf u}_{\text{v}}|^2 \right), 
$$
where $\gamma=c_p/c_v$ is the ratio of specific heat capacities at
constant pressure, $c_p$, and constant volume, $c_v$; for dry air,
$\gamma = 1.4$. Finally, $\eyedentity$ denotes the $d\times d$ identity matrix.


\begin{table}[t!]
  \begin{lstlisting}[language=Python,caption={Example 2: Specification of the convective flux and dissipation parameter $\alpha$ for the two-dimensional compressible Euler equations in UFL.}, captionpos=b,
  label=code:fc_fv_comp_euler]
def F_c(U):
  rho, u1, u2, E = U[0], U[1]/U[0], U[2]/U[0], U[3]/U[0]
  p = (gamma - 1.0)*rho*(E - 0.5*(u1**2 + u2**2))
  H = E + p/rho
  return as_matrix([[rho*u1,        rho*u2       ],
                    [rho*u1**2 + p, rho*u1*u2    ],
                    [rho*u1*u2,     rho*u2**2 + p],
                    [rho*H*u1,      rho*H*u2     ]])

def alpha(U, n):
    rho, u1, u2, E = U[0], U[1]/U[0], U[2]/U[0], U[3]/U[0]
    p = (gamma - 1.0)*rho*(E - 0.5*(u1**2 + u2**2))
    u = as_vector([u1, u2])
    c = sqrt(gamma*p/rho)
    lambdas = [dot(u, n) - c, dot(u, n), dot(u, n) + c]
    return lambdas
  \end{lstlisting}
\end{table}

The automatic construction of the DGFEM formulation for the numerical
approximation of~\eqref{eq:comp_euler_bvp} is encapsulated
within the class \ilc{CompressibleEulerOperator}. We note that this
class simply inherits the \ilc{HyperbolicOperator}
class, with the specification of the Euler flux and dissipation parameter 
required for the definition of the Lax Friedrichs flux, 
cf. Listing~\ref{code:fc_fv_comp_euler} for the case when $d=2$.
We note that in this setting 
$\lambda_{\max} = \max\{|\cvec{u}_{\text{v}}|,|\cvec{u}_{\text{v}}|\pm c\}$, 
where $c=\sqrt{\gamma p/\rho}$ is the speed of sound.

\begin{figure}[t!]
	\centering 
  \includegraphics[scale=0.3]{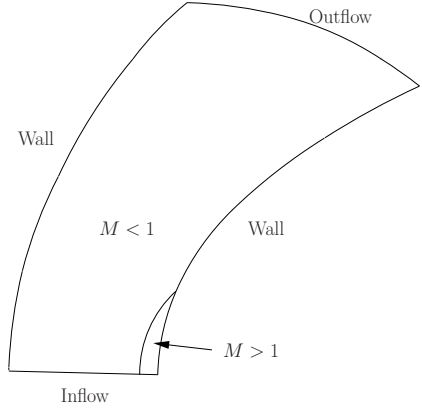}
\caption{Example 2: Computational domain for Ringleb's flow.}
\label{fig:ringleb_geom}
\end{figure}

As an illustration of the exploitation of \ilc{CompressibleEulerOperator},
we consider the DGFEM approximation of Ringleb's flow problem
for which an analytical solution may be obtained using the hodograph method,
cf.~\cite{Chiocchia1985}. The DGFEM discretisation of this test case has also been 
studied in~\cite{Hartmann2002}. This problem represents a transonic flow in a 
curved channel domain, where the flow is mainly subsonic, with a small 
supersonic region near the right-hand-side wall; cf. Figure~\ref{fig:ringleb_geom}.
Here, inflow/outflow boundary conditions are imposed on the top and bottom
boundaries of $\Omega$, while a solid wall condition is specified on the 
left- and right-hand side boundaries. Following~\cite{HH_vki_2009} the latter
boundary conditions are enforced based on employing a symmetry technique
through the specification of the boundary function ${\bf  u}_\Gamma$.
More precisely, we treat the walls as part of the Dirichlet boundary, whereby
we set
\begin{equation*}
{\bf u}_\Gamma({\bf u})=\left(
    \begin{array}{cccc}
                1 & 0 & 0 & 0\\
                0 & 1-2 n_1^2 & -2 n_1 n_2 & 0\\
                0 & -2 n_1 n_2 & 1-2 n_2^2 & 0\\
                  0 & 0 & 0 & 1           
                \end{array}
              \right){\bf u};
\end{equation*}
thereby, here $\extbdr = \extbdr[D]$ and $\extbdr[N] = \emptyset$.

\begin{figure}[t!]
\centering
\begin{subfigure}{.5\textwidth}
  \centering
  \includegraphics[width=1.\linewidth]{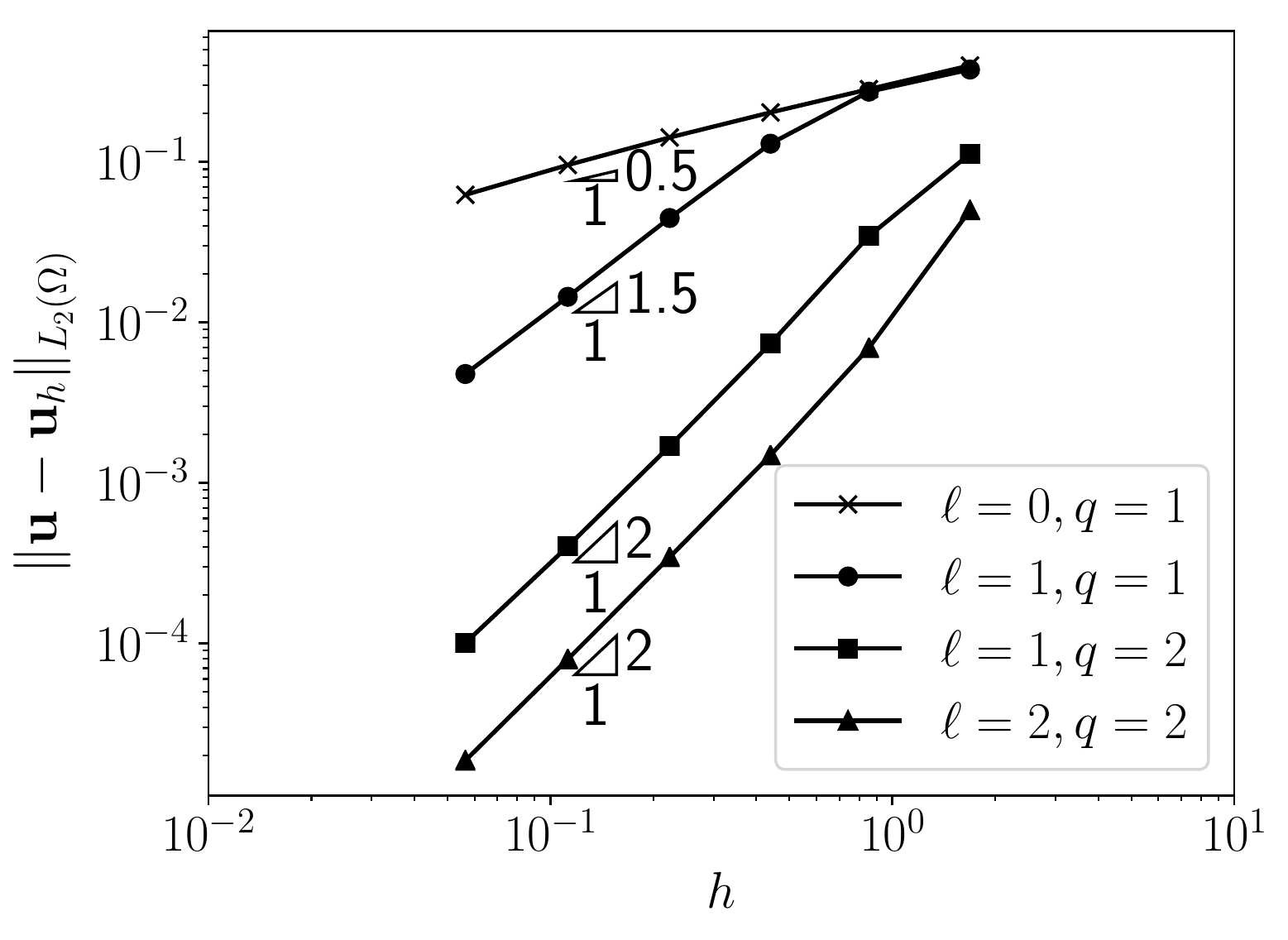} \\
  (a)
\end{subfigure}%
\begin{subfigure}{.5\textwidth}
  \centering
  \includegraphics[width=1.\linewidth]{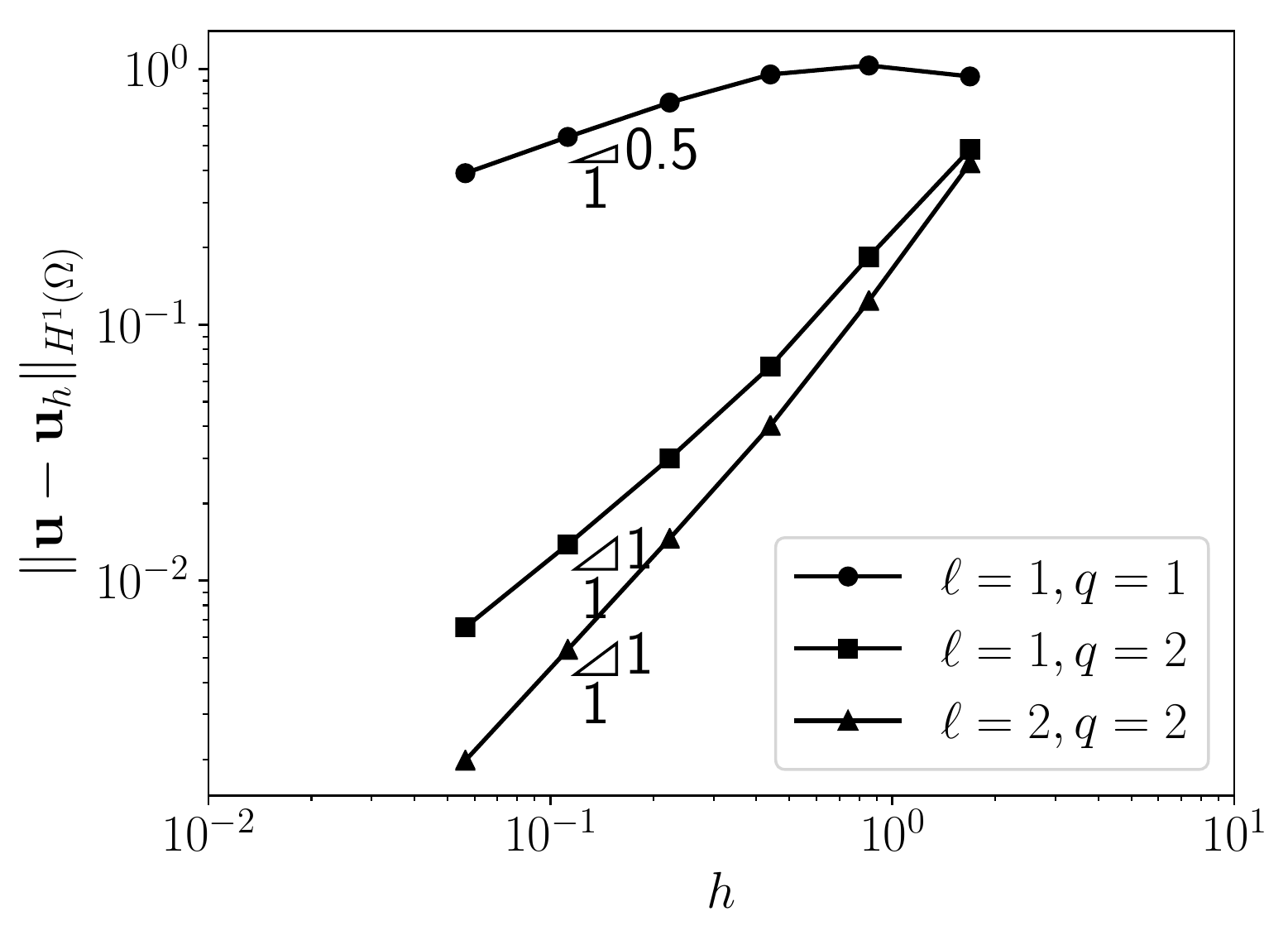} \\
  (b)
\end{subfigure}
\caption{Example 2: Convergence of the DGFEM with $h$--refinement for Ringleb's
flow: (a) $\|\cvec{u}-\cvec{u}_h\|_{L_2(\Omega)}$; 
(b) $\|\cvec{u}-\cvec{u}_h\|_{H^1(\Omega)}$.}
\label{fig:ringleb_conv_rates}
\end{figure}

\begin{table}[t!]
  \begin{lstlisting}[language=Python,caption={Example 2: Code snippet for the automatic generation of the DGFEM formulation for the numerical
  approximation of Ringleb's flow.},captionpos=b,
  label=code:comp_euler_dg_ufl]
gD = RinglebAnalyticalSoln(element=V.ufl_element(), domain=mesh)
u_vec = project(gD, V)
n = FacetNormal(mesh)

slip_proj = as_matrix(((1,            0,            0, 0),
                       (0,  1-2*n[0]**2, -2*n[0]*n[1], 0),
                       (0, -2*n[0]*n[1],  1-2*n[1]**2, 0),
                       (0,            0,            0, 1)))
slip_bc = slip_proj * u_vec

bcs = [DGDirichletBC(ds(0), gD), DGDirichletBC(ds(1), slip_bc)]
ceo = CompressibleEulerOperator(mesh, V, bcs)
residual = ceo.generate_fem_formulation(u_vec, v_vec)

solve(residual == 0, u_vec)
  \end{lstlisting}
\end{table}

Functions employed for the specification of the analytical solution
${\bf u}$ of Ringleb's flow, together with routines for the construction
of the computational mesh are provided within the file \url{ringleb.py}.
The curved boundaries on the walls of the domain must be treated
in a careful manner to ensure optimal convergence of the underlying
DGFEM discretisation. Indeed, our computational experience suggests that
a curved polynomial description of the boundary of order $q \geq \polyo + 1$ 
is necessary to ensure that the $L_2(\Omega)$
and $H^1(\Omega)$ norms of the error tend to zero at the optimal
rates of ${\mathcal O}(h^{\ell+1})$ and
${\mathcal O}(h^{\ell})$ as the mesh size $h$ tends to zero, for fixed $\polyo$,
cf. Figure~\ref{fig:ringleb_conv_rates}. The complete python code for this
numerical example is provided \url{ringleb_example.py}; a snippet of this
code is depicted in Listing~\ref{code:comp_euler_dg_ufl} in order
to highlight the simplicity of the specification of the
DGFEM for this example.


\subsection{Example 3a: Compressible Navier-Stokes Equations}
\label{subsec:comp_ns_numerical_exp}
%
In this example we consider the DGFEM discretisation of the
compressible Navier-Stokes equations:
find $\cvec{u} = (\rho,
\rho\cvec{u}_{\text{v}}, \rho E)^\top : \domain \rightarrow \mathbb{R}^+ \times
\mathbb{R}^{d} \times \mathbb{R}^+$ such that
\begin{alignat}{2}
  \nabla \cdot \left( \mathcal{F}^c \left( \cvec{u} \right) - \mathcal{F}^{v} \left(\cvec{u}, \nabla \cvec{u} \right) \right) &= \cvec{f} &&\quad\text{in } \domain, \label{eq:comp_ns_bvp} \end{alignat}
  \label{eq:comp_ns}
where $\mathcal{F}^c$ is defined as in \eqref{eq:euler_flux} and
the viscous flux is given by
\begin{equation}
  \mathcal{F}^v(\cvec{u}, \nabla \cvec{u}) = 
  \begin{pmatrix}
    \cvec{0} \\
    \tau \\
    \tau \cvec{u}_{\text{v}} + {\mathcal K} \nabla T \\
  \end{pmatrix},
\end{equation}
where $T$ denotes the temperature and ${\mathcal K}$ is the thermal conductivity
coefficient. The stress tensor is defined by
$$
\tau = \visc \left( \nabla \cvec{u}_{\text{v}} + \nabla \cvec{u}_{\text{v}}^\top - \tfrac{2}{3}
\left( \nabla \cdot \cvec{u}_{\text{v}} \right) \eyedentity \right),
$$
where $\visc$ is the dynamic viscosity coefficient.  Finally, we note that
${\mathcal K}T = \frac{\visc\gamma}{\prandtl} E - \tfrac{1}{2} \cvec{u}_{\text{v}}^2$,
where $\prandtl=0.72$ is the Prandtl number.

\begin{table}[t!]
  \begin{lstlisting}[language=Python,caption={Example 3a: Specification of the viscous flux for the two-dimensional compressible Navier-Stokes equations in UFL.},captionpos=b,
  label=code:fc_fv_comp_ns]
def F_v(U, grad_U):
    rho, rhou, rhoE = conserved_variables(U)
    u = rhou/rho

    grad_rho = grad_U[0, :]
    grad_rhou = as_matrix([[grad_U[j,:] for j in range(1, dim + 1)]])[0]
    grad_rhoE = grad_U[-1,:]
    # Quotient rule to find grad(u) and grad(E)
    grad_u = as_matrix([[(grad_rhou[j,:]*rho - rhou[j]*grad_rho)/rho**2 for j in range(dim)]])[0]
    grad_E = (grad_rhoE*rho - rhoE*grad_rho)/rho**2

    tau = mu*(grad_u + grad_u.T - 2.0/3.0*(tr(grad_u))*Identity(dim))
    K_grad_T = mu*gamma/Pr*(grad_E - dot(u, grad_u))

    r = as_matrix([[0.0,                            0.0                             ],
                   [tau[0,0],                       tau[0,1]                        ],
                   [tau[1,0],                       tau[1,1]                        ],
                   [dot(tau[0,:], u) + K_grad_T[0], (dot(tau[1,:], u)) + K_grad_T[1]]])
    return r
  \end{lstlisting}
\end{table}

\begin{figure}[t!]
\centering
\begin{subfigure}{.5\textwidth}
  \centering
  \includegraphics[width=1.\linewidth]{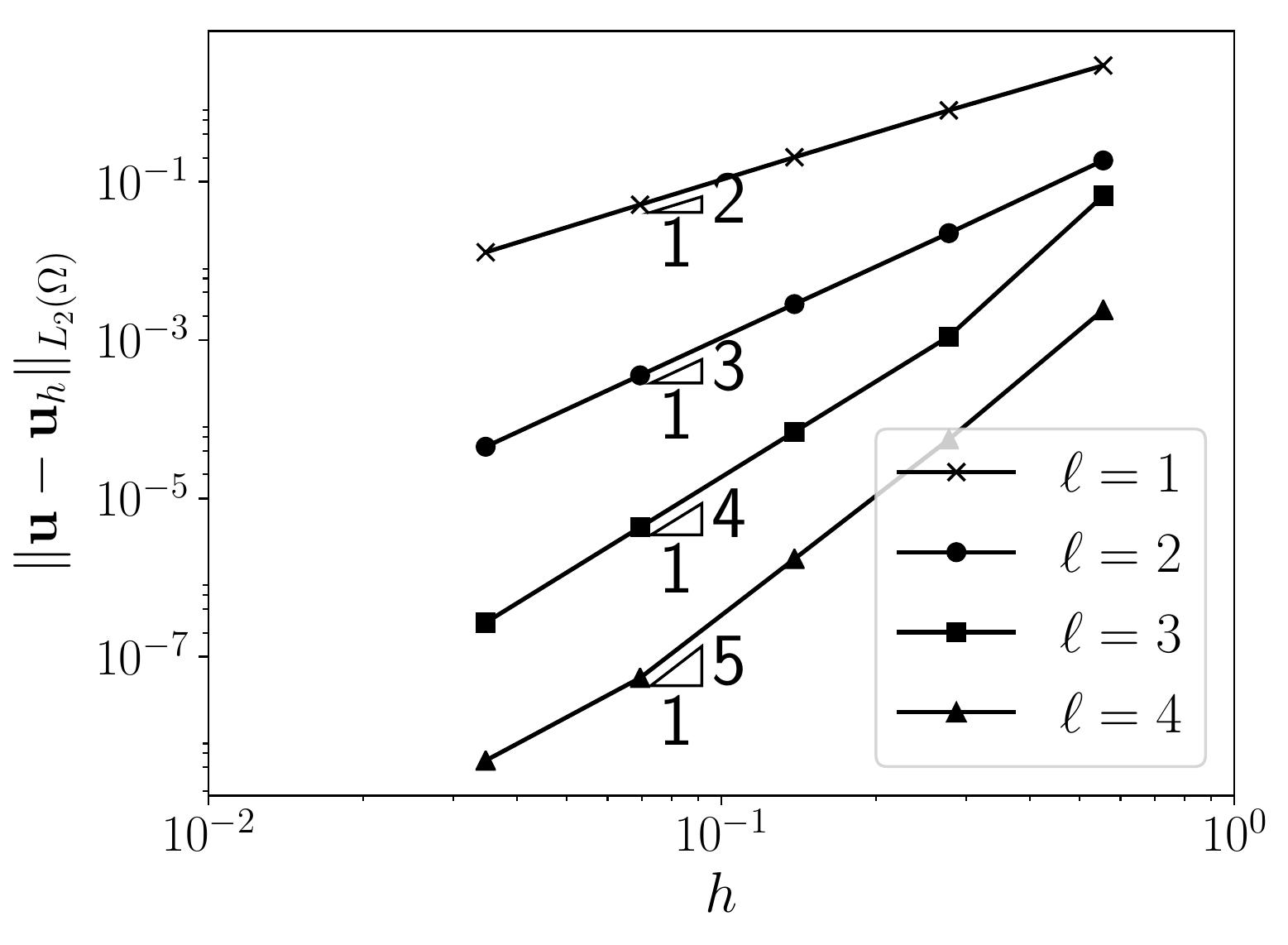}\\
  (a)
\end{subfigure}%
\begin{subfigure}{.5\textwidth}
  \centering
  \includegraphics[width=1.\linewidth]{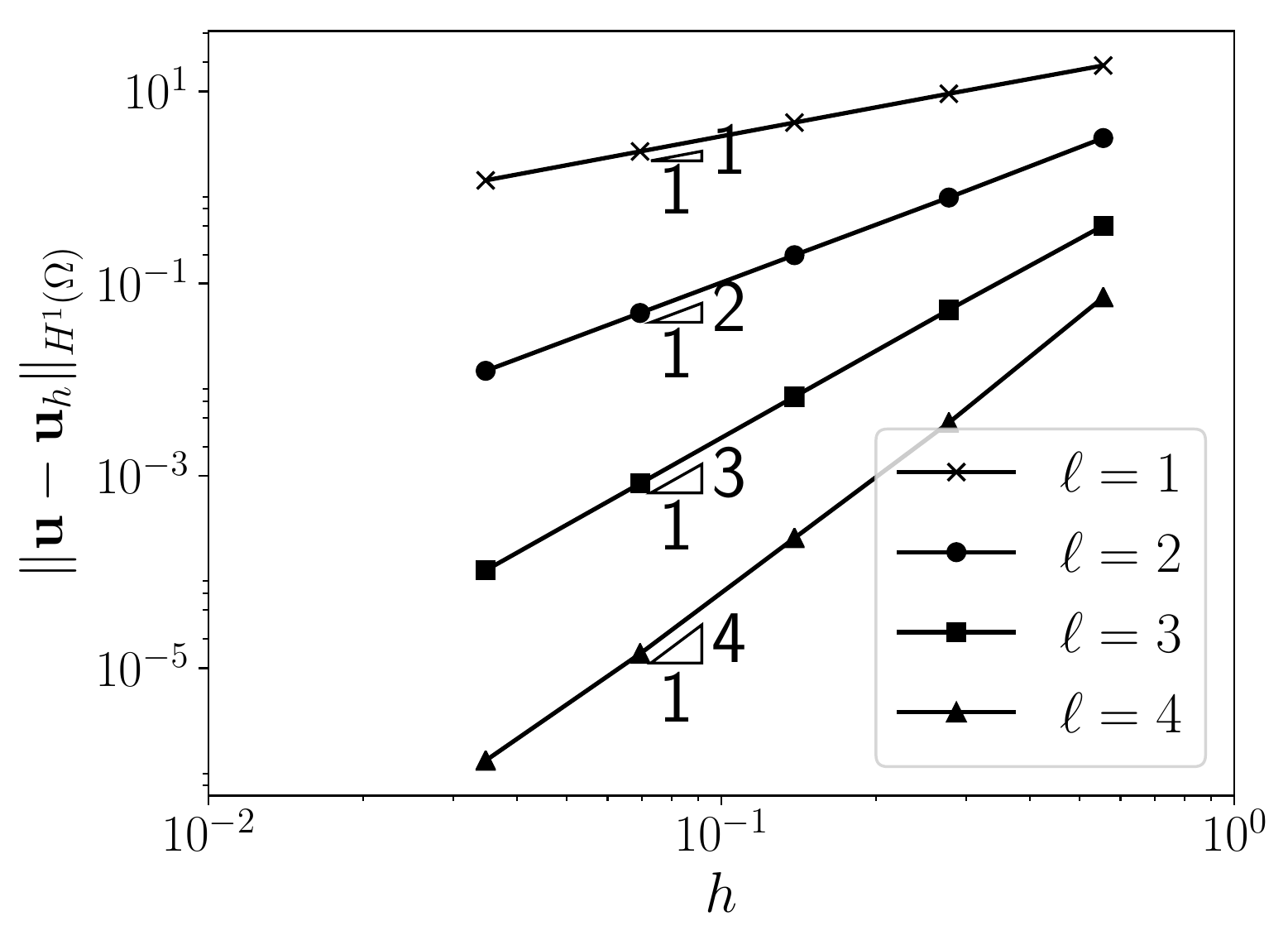}\\
  (b)
\end{subfigure}
\caption{Example 3a: Convergence of the DGFEM with $h$--refinement for the compressible Navier-Stokes equations: (a) $\|\cvec{u}-\cvec{u}_h\|_{L_2(\Omega)}$; 
(b) $\|\cvec{u}-\cvec{u}_h\|_{H^1(\Omega)}$.}
\label{fig:comp_ns_conv}
\end{figure}

Given the UFL code in Listing~\ref{code:fc_fv_comp_euler}, together with the
specification of the viscous flux in Listing~\ref{code:fc_fv_comp_ns}, in the case
when $d=2$, we have implemented the class \ilc{CompressibleNavierStokesOperator}
class, which inherits both the \ilc{CompressibleEuler\-Operator} and
\ilc{EllipticOperator} classes. Indeed, the \ilc{CompressibleEulerOperator}
component is treated in an identical manner as in the previous section, while
the \ilc{EllipticOperator} component is employed with the UFL specification of
the viscous flux. On the basis of these classes, the homogeneity tensor and the
resulting symbolic DGFEM formulation can then be automatically generated.
Moreover, the G\^{a}teaux derivative of the DGFEM formulation is automatically
computed by UFL for use within the Newton solver managed in DOLFIN by invoking
the call to \ilc{solve}. As a simple test, we consider the example outlined
in~\cite{HH08a}; namely, we set $\domain = (0, \pi)^2$ and select $\cvec{f}$ so
that the analytical solution to
\eqref{eq:comp_ns} is given by
\begin{equation} \label{eq:comp_ns_apriori_soln}
\begin{pmatrix}
\rho \\
\rho \cvec{u}_{\text{v},1} \\
\rho \cvec{u}_{\text{v},2} \\
\rho E
\end{pmatrix}
= 
\begin{pmatrix}
\sin(2(x+y)) + 4 \\
\nicefrac{1}{5}\sin(2(x+y)) + 4 \\
\nicefrac{1}{5}\sin(2(x+y)) + 4 \\
\left(\sin(2(x+y)) + 4\right)^2
\end{pmatrix},
\end{equation}
where $\cvec{u}_{\text{v}} = (\cvec{u}_{\text{v},1},\cvec{u}_{\text{v},2})^\top$.
Furthermore, we set $\visc=1$,
$\extbdr = \extbdr[D]$ and $\extbdr[N] = \emptyset$.
The snippet of UFL code required to solve this problem is given 
in Listing~\ref{code:comp_ns_dgfemviscterm}; the complete code is
provided in \url{compressible_navierstokes_square.py}. The orders of
convergence of $\|\cvec{u}-\cvec{u}_h\|_{L_2(\Omega)}$ and
$\|\cvec{u}-\cvec{u}_h\|_{H^1(\Omega)}$ are reported in Figure~\ref{fig:comp_ns_conv}
as the mesh is uniformly refined for $\polyo=1,2,3,4$; as 
in~\cite{HH08a} we observe that
$\|\cvec{u}-\cvec{u}_h\|_{L_2(\Omega)}={\mathcal O}(h^{\polyo+1})$ and
$\|\cvec{u}-\cvec{u}_h\|_{H^1(\Omega)}={\mathcal O}(h^{\polyo})$
as $h$ tends to zero, for each fixed polynomial order $\polyo$.
\begin{table}[t!]
  \begin{lstlisting}[language=Python,caption={Example 3a: Code snippet for the automatic generation of the DGFEM formulation for the numerical approximation of the
  compressible Navier-Stokes equations.},captionpos=b,label=code:comp_ns_dgfemviscterm]
    gD = Expression(('sin(2*(x[0]+x[1])) + 4',
                     '0.2*sin(2*(x[0]+x[1])) + 4',
                     '0.2*sin(2*(x[0]+x[1])) + 4',
                     'pow((sin(2*(x[0]+x[1])) + 4), 2)'),
                    element=V.ufl_element())

    f = Expression((...),
                   element=V.ufl_element())

    u, v = interpolate(gD, V), TestFunction(V)

    cnso = CompressibleNavierStokesOperator(mesh, V, DGDirichletBC(ds, gD))
    residual = cnso.generate_fem_formulation(u, v) - inner(f, v)*dx

    solve(residual == 0, u)
  \end{lstlisting}
\end{table}

\subsection{Example 3b: Compressible Flow Around a NACA0012 Airfoil}
\label{subsec:comp_ns_drag_numerical_exp}

To demonstrate the application of the automatic formulation of DGFEMs 
applied to exterior flow problems, in this section
we consider laminar flow around a NACA0012 airfoil whose
geometry is defined by
\begin{equation}
y = \pm 5t\left(0.2969\sqrt{x} - 0.1260 x - 0.3516 x^2 + 0.2843 x^3 - 0.1036 x^4\right),
\end{equation}
where the thickness fraction $t = \num{0.12}$.
Here, we set the angle of attack $\alpha = 2^\circ$,
Reynolds number $\mathrm{Re} = \num{5e3}$, inlet flow Mach number $M=0.5$, 
and impose an adiabatic no slip condition on the airfoil; 
cf.~\cite{HH08a,Hartmann2005}. We note that the adiabatic no slip
condition imposes a zero heat
flux condition $\mathcal{K} \nabla T \cdot \mathbf{n} = 0$ on the airfoil
surface boundary. We have implemented this by 
defining a new boundary condition in the
class \ilc{DGAdiabaticWallBC} inheriting \ilc{DGBC}. This boundary condition is
then recognised by the \ilc{CompressibleNavierStokesOperator} implementation
which automatically generates the DGFEM formulation accordingly.

\begin{figure}[t!]
\centering
\begin{subfigure}{.5\textwidth}
  \centering
  \includegraphics[width=1.\linewidth]{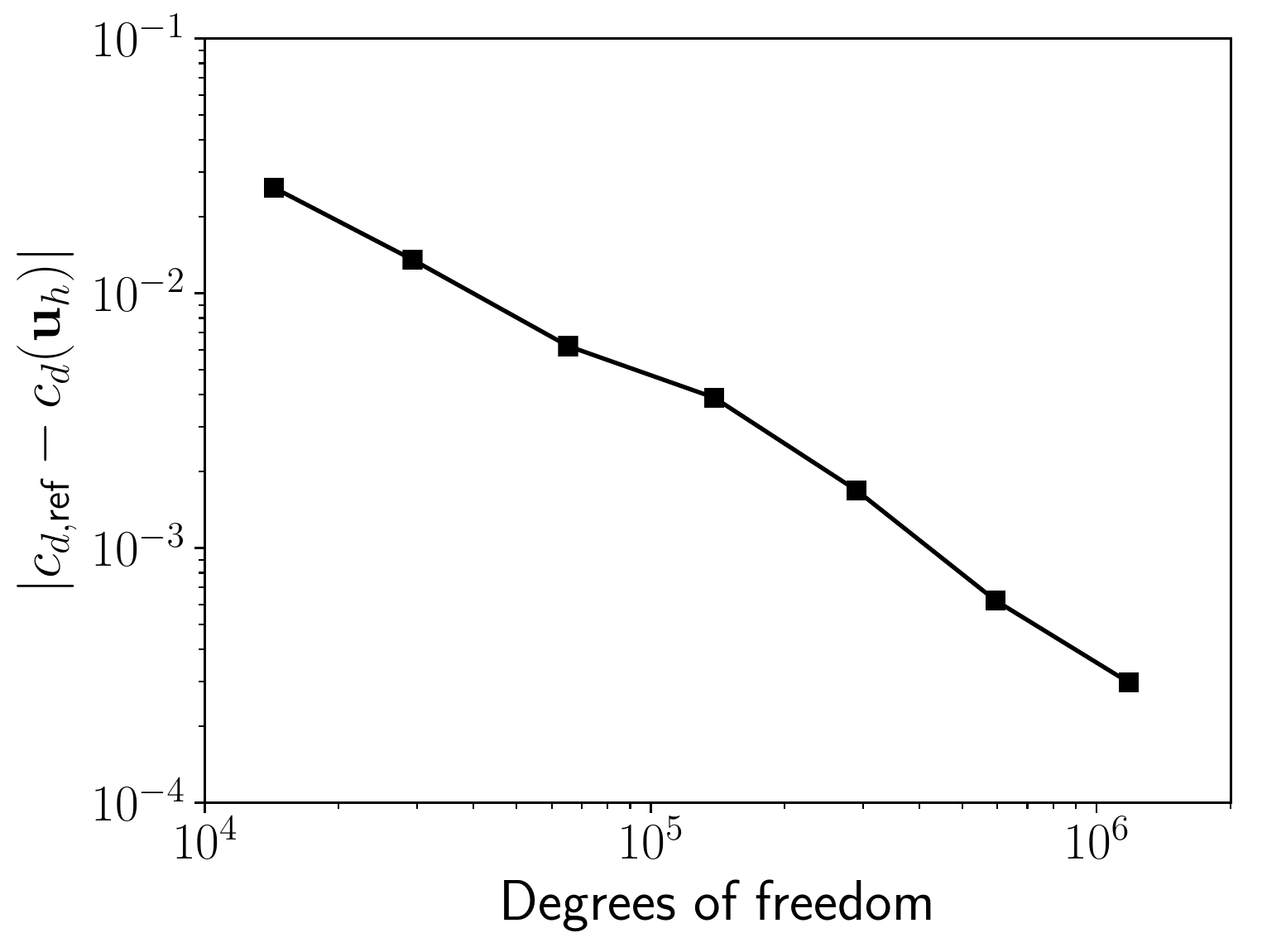}\\
  (a)
\end{subfigure}
\begin{subfigure}{.5\textwidth}
  \centering
  \includegraphics[width=1.\linewidth]{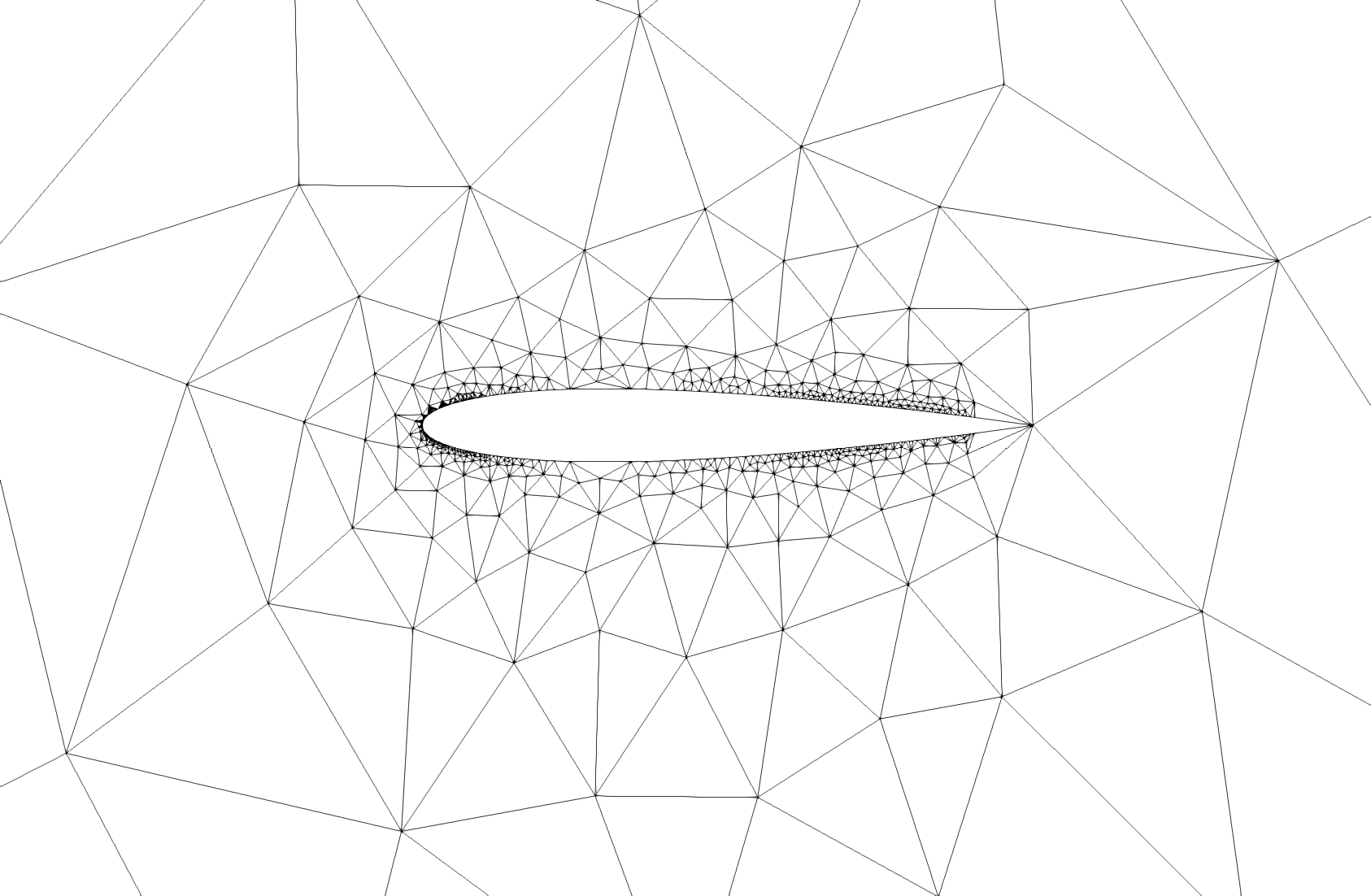}\\
  (b)
\end{subfigure}%
\begin{subfigure}{.5\textwidth}
  \centering
  \includegraphics[width=1.\linewidth]{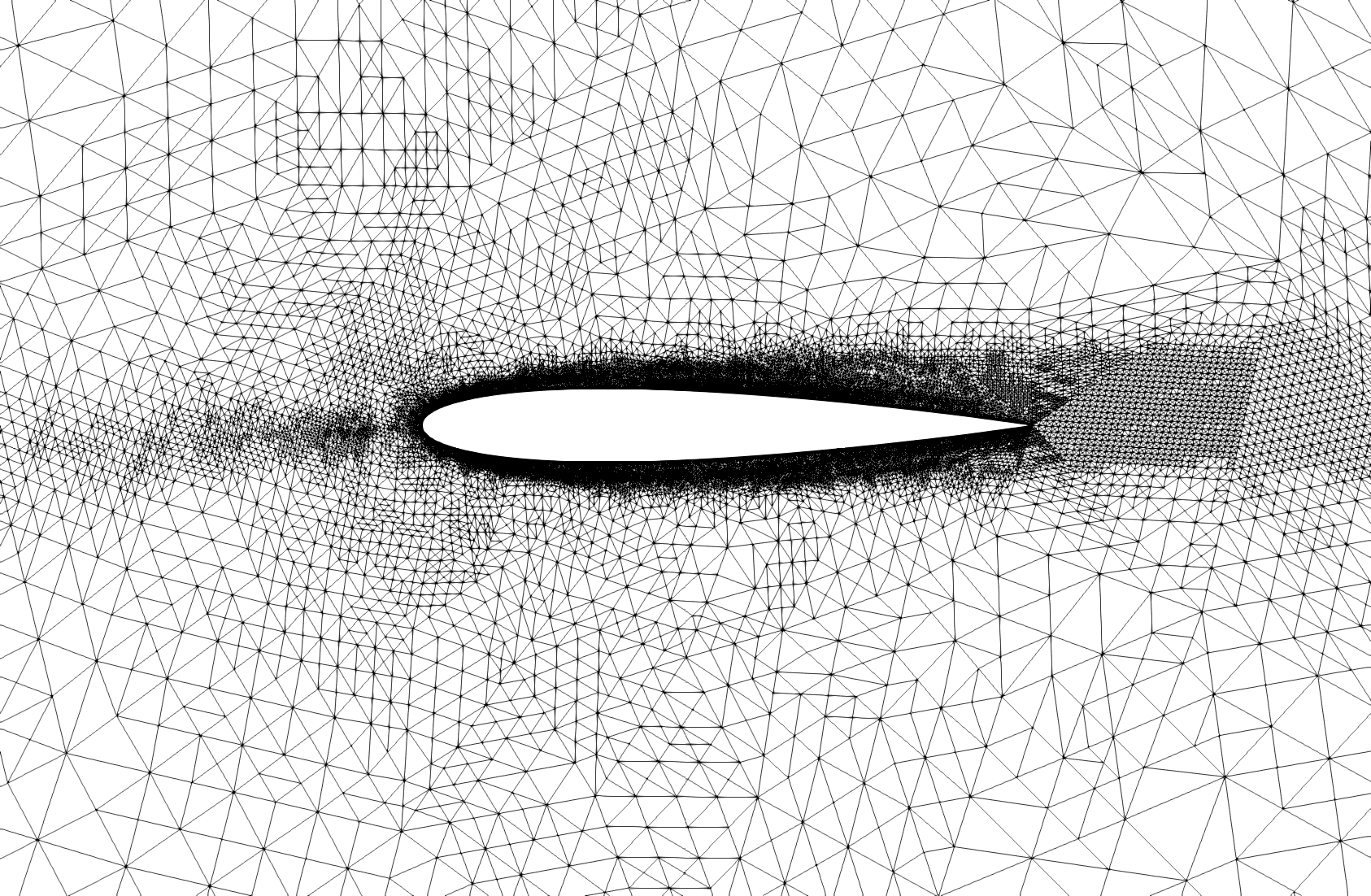}\\
  (c)
\end{subfigure}
\caption{Example 3b: Convergence of the DGFEM with DWR $h$--refinement for the
compressible Navier-Stokes equations with $\ell=1$: (a) $\left|C_{d,\text{ref}}-C_{d}(\cvec{u}_h)\right|$; (b) Initial mesh; (c) Final mesh.}
\label{fig:comp_ns_drag_conv}
\end{figure}

In this example, we shall undertake goal-oriented adaptive mesh
refinement based on employing a dual
weighted residual (DWR) \emph{a posteriori} error estimator;
for details, we refer to~\cite{Hartmann2005}. In particular, we select
the quantity of interest to be the drag coefficient
\begin{equation}
c_d(\mathbf{u}) = \int_{\Gamma_{\text{airfoil}}} \frac{1}{C_\infty} \left( p \mathbf{n}  - \tau \mathbf{n} \right) \cdot \psi_\text{drag} \; \mathrm{d}s,
\end{equation}
where $\Gamma_{\text{airfoil}}$ denotes the surface of the airfoil,
$\psi_\text{drag} = (\cos(\alpha), \sin(\alpha))^\top$, $C_\infty = 
\frac{1}{2} \rho_{\infty} \mathbf{u}_{\mathbf{v},\infty}^2 l_\text{ref}$,
$\rho_{\infty}$ and $\mathbf{u}_{\mathbf{v},\infty}$ denote the
far field density and velocity, respectively, and $l_\text{ref}$ is the
reference length scale. The symbolic representation of the
underlying dual problem is automatically formulated using
the \texttt{dwr} component of \texttt{dolfin\_dg}. The error 
in the approximate drag coefficient on each of the meshes employed 
is evaluated relative to a reference drag
coefficient, $c_{d,\text{ref}}$ computed from a very fine mesh problem. Setting
$\ell=1$, the convergence of the error in the computed drag coefficient with
respect to the number of degrees of freedom in the underlying finite element
space, together with the initial and final computational meshes
are shown in Figure~\ref{fig:comp_ns_drag_conv}. The code to generate these
results is available in the file 
\texttt{dg\_naca0012\_2d.py}. 

\subsection{Example 3c: Entropy Variable Formulation of the Compressible Navier-Stokes Equations}

\begin{table}[t!]
  \begin{lstlisting}[language=Python,caption={Example 3c: Transformations~\eqref{eq:u_to_v} and~\eqref{eq:v_to_u} in UFL representation.}, captionpos=b, label=code:comp_ns_entropy_mapping]{}
def U_to_V(U, gamma):
    rho, u1, u2, E = U[0], U[1]/U[0], U[2]/U[0], U[3]/U[0]
    i = E - 0.5*(u1**2 + u2**2)
    U1, U2, U3, U4 = U
    s = ln((gamma-1)*rho*i/(U1**gamma))
    V1 = 1/(rho*i)*(-U4 + rho*i*(gamma + 1 - s))
    V2 = 1/(rho*i)*U2
    V3 = 1/(rho*i)*U3
    V4 = 1/(rho*i)*(-U1)
    return as_vector([V1, V2, V3, V4])

def V_to_U(V, gamma):
    V1, V2, V3, V4 = V
    U = as_vector([-V4, V2, V3, 1 - 0.5*(V2**2 + V3**2)/V4])
    s = gamma - V1 + (V2**2 + V3**2)/(2*V4)
    rhoi = ((gamma - 1)/((-V4)**gamma))**(1.0/(gamma-1))*exp(-s/(gamma-1))
    U = U*rhoi
    return U
  \end{lstlisting}
\end{table}

\begin{table}[t!]
  \begin{lstlisting}[language=Python,caption={Example 3c: Convective and viscous fluxes of the entropy formulation of the compressible Navier-Stokes equations. For the sake of brevity we do not
  repeat the code for the viscous flux tensor which is provided in
  Listing~\ref{code:fc_fv_comp_ns}.}, captionpos=b,
  label=code:fc_fv_comp_ns_entropy]
def F_c(V):
    U = V_to_U(V, gamma)
    rho, u1, u2, E = U[0], U[1]/U[0], U[2]/U[0], U[3]/U[0]
    p = (gamma - 1.0)*rho*(E - 0.5*(u1**2 + u2**2))
    H = E + p/rho
    res = as_matrix([[rho*u1,        rho*u2       ],
                     [rho*u1**2 + p, rho*u1*u2    ],
                     [rho*u1*u2,     rho*u2**2 + p],
                     [rho*H*u1,      rho*H*u2     ]])
    return res

def F_v(V, grad_V):
    V = variable(V)
    U = V_to_U(V, gamma)
    dudv = diff(U, V)
    grad_U = dot(dudv, grad_V)

    ...

    r = as_matrix([[0.0,0.0],
                   [tau[0,0],tau[0,1]],
                   [tau[1,0],tau[1,1]],
                   [dot(tau[0,:], u) + K_grad_T[0], (dot(tau[1,:], u)) + K_grad_T[1]]])
    return r
  \end{lstlisting}
\end{table}

\begin{table}[t!]
  \begin{lstlisting}[language=Python,caption={Example 3c: Code snippet for the automatic generation of the DGFEM
  formulation for the numerical approximation of the entropy formulation of the
  compressible Navier-Stokes equations.}, captionpos=b,
  label=code:comp_ns_entropy_dgfemviscterm]
    # Dirichlet conditions and error suite
    gamma = 1.4
    gD = Expression((...),
                    element=V.ufl_element())

    f = Expression((...),
                   element=V.ufl_element())

    V_vec, V_test = interpolate(gD, V), TestFunction(V)

    bo = CompressibleNavierStokesOperatorEntropyFormulation(mesh, V, DGDirichletBC(ds, gD))
    residual = bo.generate_fem_formulation(V_vec, V_test) - inner(f, V_test)*dx

    solve(residual == 0, V_vec)
  \end{lstlisting}
\end{table}

In this section, we consider the DGFEM approximation of the compressible
Navier-Stokes equations written in terms of so--called (symmetrisation) entropy
variables; for further details, we refer, for example,
to~\cite{Barth_1999,Hughes1986}. For simplicity of presentation we only consider
the two--dimensional case, though the extension to $d=3$ follows in an analogous
manner. Thereby, we introduce the change of variable $\cvec{u}
\mapsto \cvec{V}(\cvec{u})$ given by
\begin{equation} \label{eq:u_to_v}
\cvec{V} = \frac{1}{\rho e}
\begin{pmatrix}
-\rho E + \rho e ( \gamma + 1 - s) \\
\rho \cvec{u}_{\text{v},1} \\
\rho \cvec{u}_{\text{v},2} \\
-\rho
\end{pmatrix},
\;
s = \ln \left( \frac{(\gamma - 1) \rho e}{\rho^\gamma} \right),
\;
e = E - \frac{1}{2} |\cvec{u}_{\text{v}}|^2,
\end{equation}
and its inverse
\begin{eqnarray} 
\cvec{u} &=& \rho e
\begin{pmatrix}
-V_4 \\
V_2 \\
V_3 \\
1 - \frac{1}{2}(V_2^2 + V_3^2)/V_4
\end{pmatrix},
\;
\rho e = \left( \frac{\gamma - 1}{(-V_4)^\gamma} \right)^\frac{1}{\gamma-1} \exp \left( \frac{-s}{\gamma - 1} \right), 
\nonumber \\
s &=& \gamma - V_1 + \frac{1}{2}\frac{(V_2^2 + V_3^2)}{V_4}.\label{eq:v_to_u}
\end{eqnarray}
Here, $e$ is the internal energy density and $s$ is a nondimensional entropy. 
The implementation of this mapping and its inverse in 
the UFL symbolic algebra framework is depicted in 
Listing~\ref{code:comp_ns_entropy_mapping}. On the basis of
this transformation, the symmetrised formulation of the
compressible Navier-Stokes equations is given by: find
$\cvec{V} = (V_1, V_2, V_3, V_4)^\top : \domain \rightarrow \mathbb{R}^- \times
\mathbb{R}^{2} \times \mathbb{R}^-$ such that
\begin{alignat}{2} \label{eq:comp_ns_entropy_bvp}
\nabla \cdot \left( \mathcal{F}^c_{\cvec{V}} ( \cvec{V} ) - \mathcal{F}^v_{\cvec{V}}( \cvec{V}, \nabla \cvec{V}) \right) &= \cvec{f} &&\quad\text{in } \domain.\end{alignat}
where $\mathcal{F}^c_{\cvec{V}} ( \cvec{V} ) = \mathcal{F}^c ( \cvec{u}(\cvec{V}) )$
and $\mathcal{F}^v_{\cvec{V}}( \cvec{V}, \nabla \cvec{V}) 
= \mathcal{F}^v( \cvec{u}(\cvec{V}), \nabla \cvec{u}(\cvec{V}))$.

\begin{figure}[t!]
\centering
\begin{subfigure}{.5\textwidth}
  \centering
  \includegraphics[width=1.\linewidth]{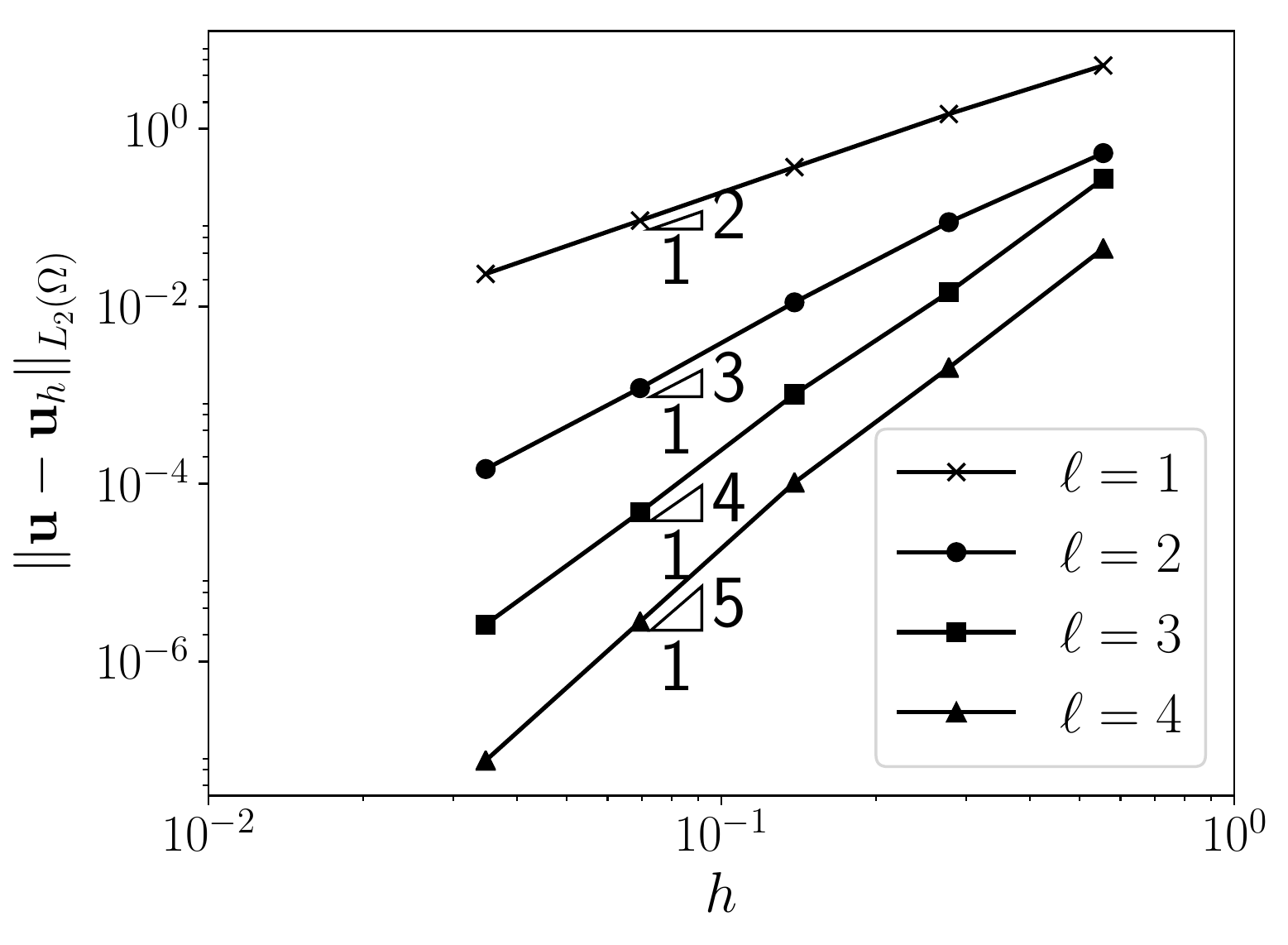}\\
  (a)
\end{subfigure}%
\begin{subfigure}{.5\textwidth}
  \centering
  \includegraphics[width=1.\linewidth]{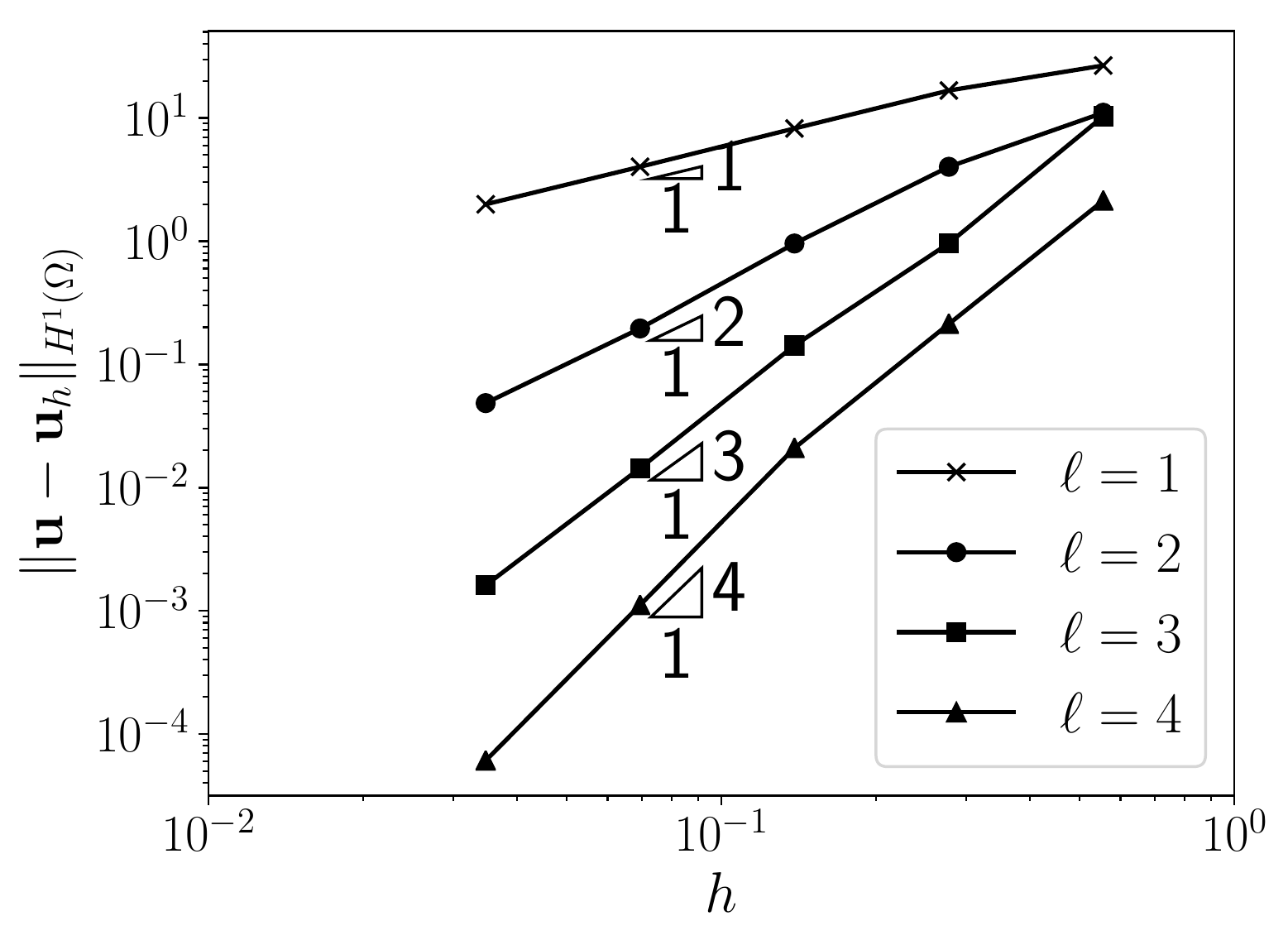}\\
  (b)
\end{subfigure}
\caption{Example 3c: Convergence of the DGFEM with $h$--refinement for the compressible Navier-Stokes equations in entropy formulation: (a) $\|\cvec{u}-\cvec{u}_h\|_{L_2(\Omega)}$; 
(b) $\|\cvec{u}-\cvec{u}_h\|_{H^1(\Omega)}$.}
\label{fig:comp_ns_entropy_conv}
\end{figure}

On the basis of the compressible Euler and compressible Navier-Stokes examples
presented in the previous two sections, the DGFEM formulation of
problem~\eqref{eq:comp_ns_entropy_bvp} is straightforward in the symbolic
algebra framework of the UFL. We use the convective and viscous fluxes shown in
Listings~\ref{code:fc_fv_comp_euler} and~\ref{code:fc_fv_comp_ns} and implement
the transformations~\eqref{eq:u_to_v} and~\eqref{eq:v_to_u} noting that $\nabla
\mathbf{u}(\mathbf{V}) = \frac{\partial \mathbf{u}(\mathbf{V})}{\partial
\mathbf{V}} \nabla \mathbf{V}$. Indeed, the UFL formulation of the convective
and viscous fluxes, $\mathcal{F}^c_{\cvec{V}}$ and $\mathcal{F}^v_{\cvec{V}}$,
respectively, are presented in Listing~\ref{code:fc_fv_comp_ns_entropy}. These
constructs are then exploited within the DGFEM utility functions in the same
manner as in the examples presented in Listings~\ref{code:comp_euler_dg_ufl}
and~\ref{code:comp_ns_dgfemviscterm}. As a simple test, we consider the
numerical approximation of the compressible Navier-Stokes example presented in
Section~\ref{subsec:comp_ns_numerical_exp}, based on employing the above entropy
variable formulation. To this end, the underlying test problem is first
transformed according to the change of variable $\cvec{u} \mapsto
\cvec{V}(\cvec{u})$, cf. ~\eqref{eq:u_to_v}, whereby the DGFEM is computed; we
then transform the numerical solution according to the inverse
mapping~\eqref{eq:v_to_u} in order to evaluate the error in the underlying
approximation. As in the previous example, we again observe optimal convergence
of the underlying DGFEM, cf. Figure~\ref{fig:comp_ns_entropy_conv}; the code
for this example is provided in
Listing~\ref{code:comp_ns_entropy_dgfemviscterm}, 
cf. \url{compressible_navierstokes_entropy_square.py}.

\subsection{Example 4: Maxwell Operator} \label{sec:maxwell}

\begin{figure}[t!]
\centering
\begin{subfigure}{.5\textwidth}
  \centering
  \includegraphics[width=1.\linewidth]{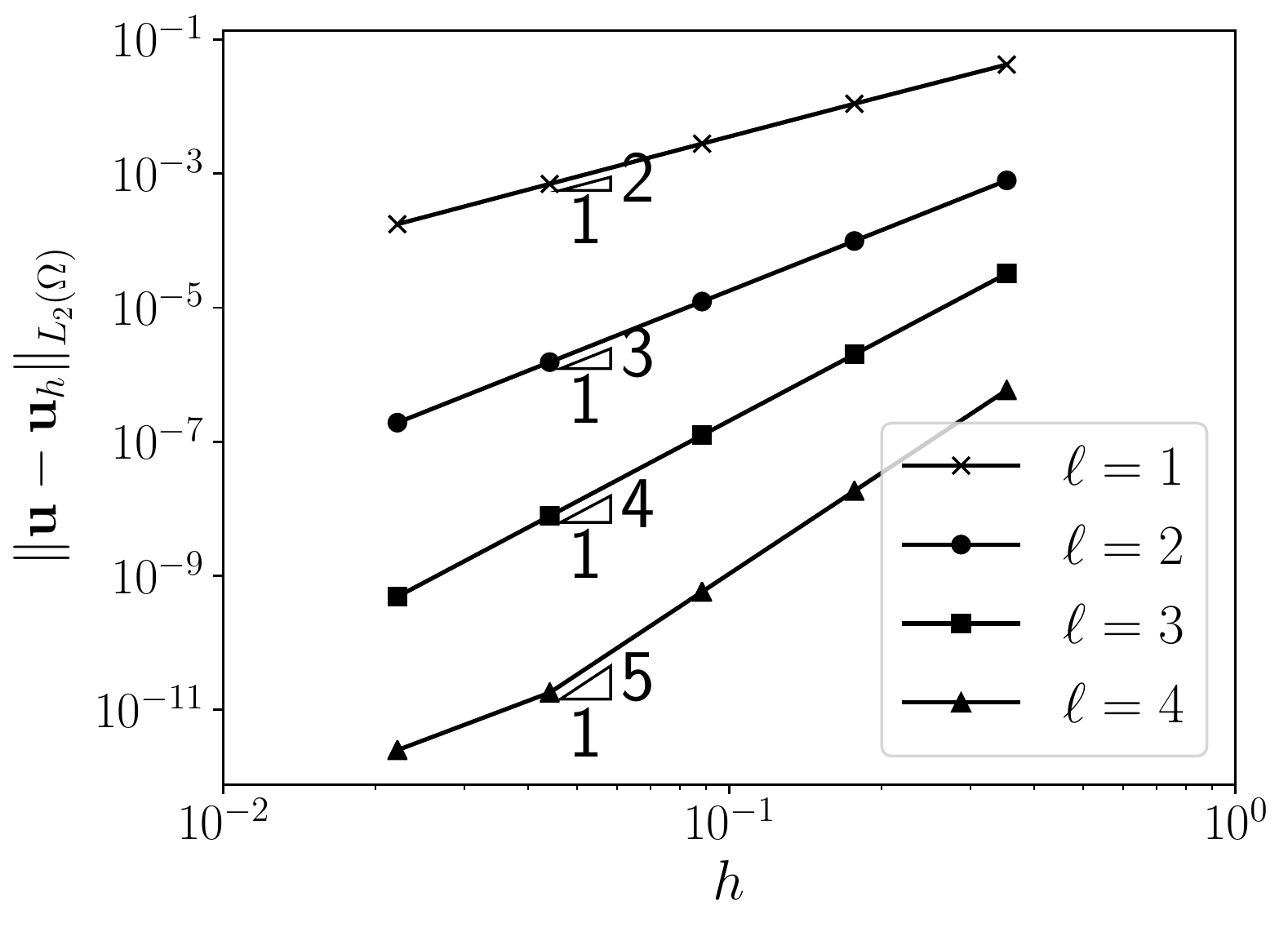} \\
  (a)
\end{subfigure}%
\begin{subfigure}{.5\textwidth}
  \centering
  \includegraphics[width=1.\linewidth]{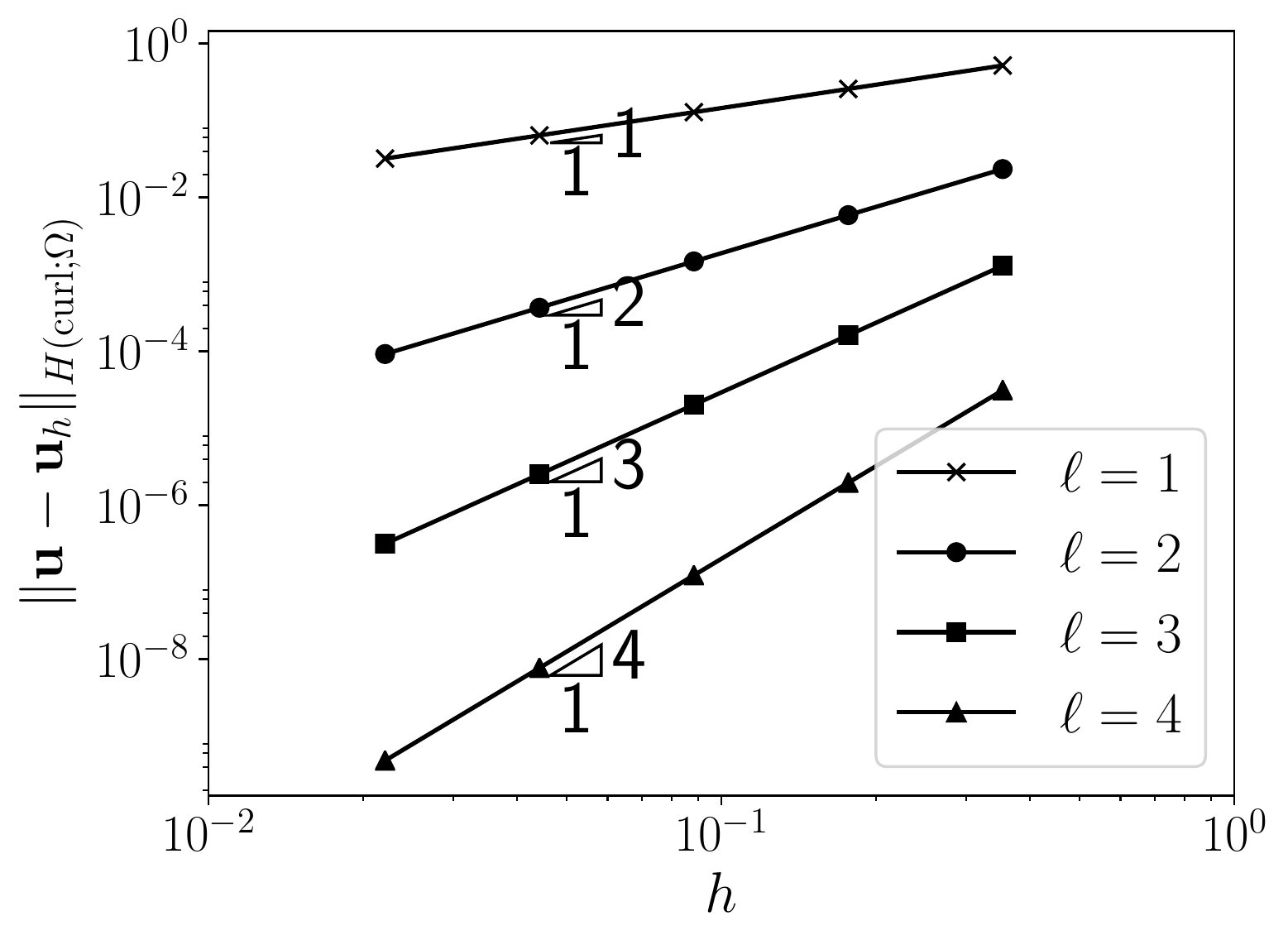} \\
  (b)
\end{subfigure}
\caption{Example 4: Convergence of the DGFEM with $h$--refinement: 
(a) $\|\cvec{u}-\cvec{u}_h\|_{L_2(\Omega)}$; (b) $\|\cvec{u}-\cvec{u}_h\|_{H(\mathrm{curl}; \Omega)}$.}
\label{fig:ex5_ad_conv_rates}
\end{figure}

In this penultimate example we demonstrate the extensibility of the
symbolic DGFEM framework developed in this article by considering the
automatic discretisation of the Maxwell problem: find $\cvec{u}$ such
that
\begin{alignat}{2}
  \nabla \times \mathcal{F}^{m} \left(\cvec{u}, \nabla \times \cvec{u} \right) - k^2 \cvec{u} &= \cvec{f} &&\quad\text{in } \domain,  \label{eq:maxwell} \\
  \cvec{n} \times \cvec{u} &= \cvec{n} \times \cvec{g}_D && \quad\text{on } \Gamma, \label{eq:maxwell_bc}
\end{alignat}
where, for simplicity, we set $\mathcal{F}^m \left(\cvec{u}, \nabla \times
\cvec{u} \right) = \nabla \times \cvec{u}$ and $k>0$ is the wave number; for
solvability, we assume that $k^2$ is not a Maxwell eigenvalue. The class
\ilc{MaxwellOperator} implements the symmetric interior penalty
discretisation of \eqref{eq:maxwell}, \eqref{eq:maxwell_bc} outlined
in~\cite{HoustonPerugiaSchneebeliSchoetzau_2005}. Thereby, given
%
%
the ufl representation of $\mathcal{F}^m \left(\cvec{u}, \nabla \times \cvec{u}
\right)$ depicted in Listing~\ref{code:ufl_maxwell_fm}, the corresponding code 
needed to compute the DGFEM approximation of  \eqref{eq:maxwell},
\eqref{eq:maxwell_bc} can be written in the very compact form shown in
Listing~\ref{code:ufl_maxwell_full}; see \url{maxwell.py}. As a simple
test case, here we have set $\domain = (-1, 1)^2$, $k=2$, $\cvec{g}_D = (\sin(k
x), \sin(k y))^\top$ and $\cvec{f} = \cvec{0}$ such that the analytical solution
is given by $\cvec{u} = \cvec{g}_D$. Numerical experiments demonstrating the
performance of the resulting scheme on a sequence of uniformly refined
triangular meshes are presented in Figure~\ref{fig:ex5_ad_conv_rates}.

\begin{table}[t!]
\begin{lstlisting}[language=Python,caption={Example 4: UFL representation of $\mathcal{F}^m\left(\cvec{u}, \nabla \times \cvec{u} \right)$ of~\eqref{eq:maxwell}.}, captionpos=b, label=code:ufl_maxwell_fm]
def F_m(u, curl_u):
    return curl_u
\end{lstlisting}
\end{table}
\begin{table}[t!]{}
\begin{lstlisting}[language=Python, caption={Example 4: Automatic DGFEM formulation for the numerical approximation of~\eqref{eq:maxwell}, \eqref{eq:maxwell_bc}.}, captionpos=b, label=code:ufl_maxwell_full]
mesh = RectangleMesh(Point(-1., -1.), Point(1., 1.), 32, 32)

V = VectorFunctionSpace(mesh, "DG", 1)
u, v = Function(V), TestFunction(V)

k = Constant(2.0)
gD = Expression(("sin(k*x[1])", "sin(k*x[0])"), k=k, element=V.ufl_element())

mo = MaxwellOperator(mesh, V, [DGDirichletBC(ds, gD)], F_m)
residual = mo.generate_fem_formulation(u, v) - k**2*dot(u, v)*dx

solve(residual == 0, u)
\end{lstlisting}
\end{table}

\subsection{Example 5a: Hyperelasticity} \label{sec:hyperelasticity}

Our final two examples highlight the flexibility of the DGFEM framework proposed
in this article for the discretisation of hyperelasticity problems; in this setting
DGFEM schemes offer computational benefits in the nearly--incompressible 
regime, cf.~\cite{HANSBO20021895,Eyck2006}. Given a domain
$\Omega_0$ defining an elastic body's reference configuration with boundary
$\partial \Omega_0$ and outward pointing unit normal vector $\mathbf{N}$, we
seek the displacement vector at all points in the domain $\mathbf{u}(\mathbf{X})
: \Omega_0 \rightarrow \mathbb{R}^3$ which determines the mapping from the
reference $\mathbf{X} \in \Omega_0$ to the deformed $\mathbf{x} \in \Omega$
configuration, such that $\mathbf{x} = \mathbf{X} + \mathbf{u}$. The
constitutive model demands mass balance, so we seek $\mathbf{u}$ such that
\begin{subequations}
\label{eq:hyperelasticity}
\begin{align}
- \nabla_{\mathbf{X}} \cdot \mathbf{P} &= \mathbf{B} \quad \text{in} \; \Omega_0, \\
\mathbf{u} &= \mathbf{u_0} \quad \text{on} \; \partial \Omega_{0,D}, \\
\mathbf{P} \cdot \mathbf{N} &= \mathbf{T} \quad \text{on} \; \partial \Omega_{0,N},
\end{align}
\label{eq:hyperelasticity_piola_kirchoff}
\end{subequations}
where $(\nabla_\mathbf{X})_i \defeq \partial / \partial X_i$ is
the gradient in the reference domain, $\mathbf{P} = \partial
\Psi(\mathbf{F}) / \partial \mathbf{F}$ is the first Piola-Kirchoff stress
tensor, $\Psi(\cdot)$ is the strain-energy function, $\mathbf{B}$ is a body
force, $\mathbf{F} = \mathbf{I} + \nabla_{\mathbf{X}} \mathbf{u}$ is the 
strain tensor, $\mathbf{T}$ is a traction force, and 
$\partial \Omega_{0,D}$ and $\partial \Omega_{0,N}$ denote the Dirichlet and
Neumann boundaries, respectively. Here,
we examine a neo-Hookean hyperelasticity model where the strain energy
density function is given by
\begin{equation*}
\Psi(\mathbf{F}) = 
\frac{\mu}{2} \left(\mathrm{tr}(\mathbf{F}^\top \mathbf{F}) - 3\right) 
- \mu \ln (\det \mathbf{F})
+ \frac{\lambda}{2} \left( \ln (\det \mathbf{F}) \right)^2,
\end{equation*}
where $\mu = E / \left(2\left( 1 + \nu \right)\right)$ and $\lambda = E \nu /
\left(\left(1 + \nu\right) \left(1 - 2 \nu \right)\right)$ are the first and
second Lam\'{e} parameters, respectively, given in terms of Young's modulus $E$
and Poisson ratio $\nu$.

\begin{table}[t!]
\begin{lstlisting}[language=Python,caption={Example 5a: UFL representation of 
the neo-Hookean hyperelasticity
equation~\eqref{eq:hyperelasticity_piola_kirchoff}.}, captionpos=b,
label=code:ufl_neohookean]
def F_v(u, grad_u):
    F = variable(Identity(3) + grad_u)
    psi = (mu/2)*(tr(F.T*F) - 3) - mu*ln(det(F)) + (lmbda/2)*(ln(det(F)))**2
    return diff(psi, F)
\end{lstlisting}
\end{table}

\begin{figure}[t!]
\centering
\begin{subfigure}{.5\textwidth}
\centering
  \includegraphics[width=1.\linewidth]{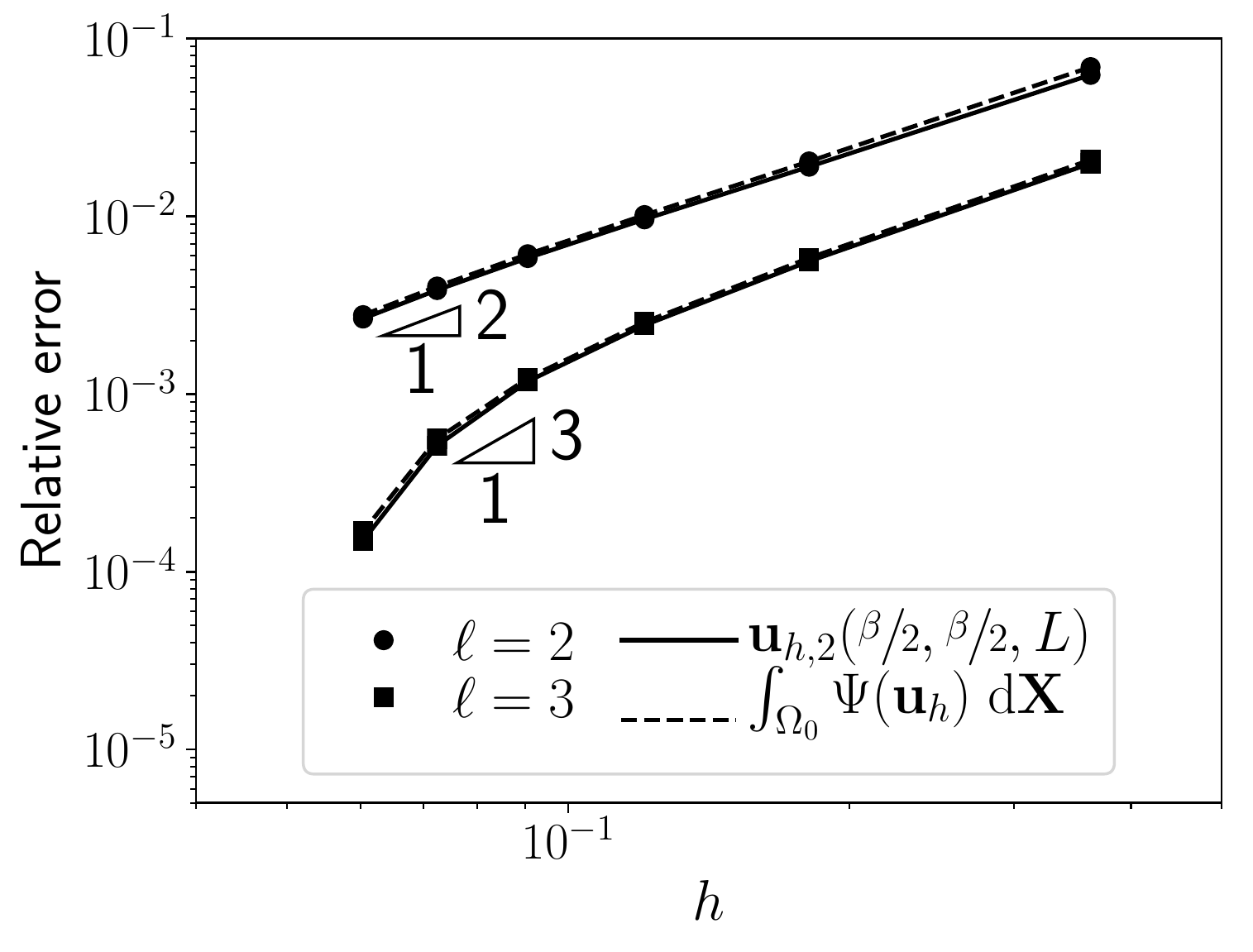}\\
  (a)
\end{subfigure}%
\begin{subfigure}{.5\textwidth}
\centering
  \includegraphics[width=1.\linewidth]{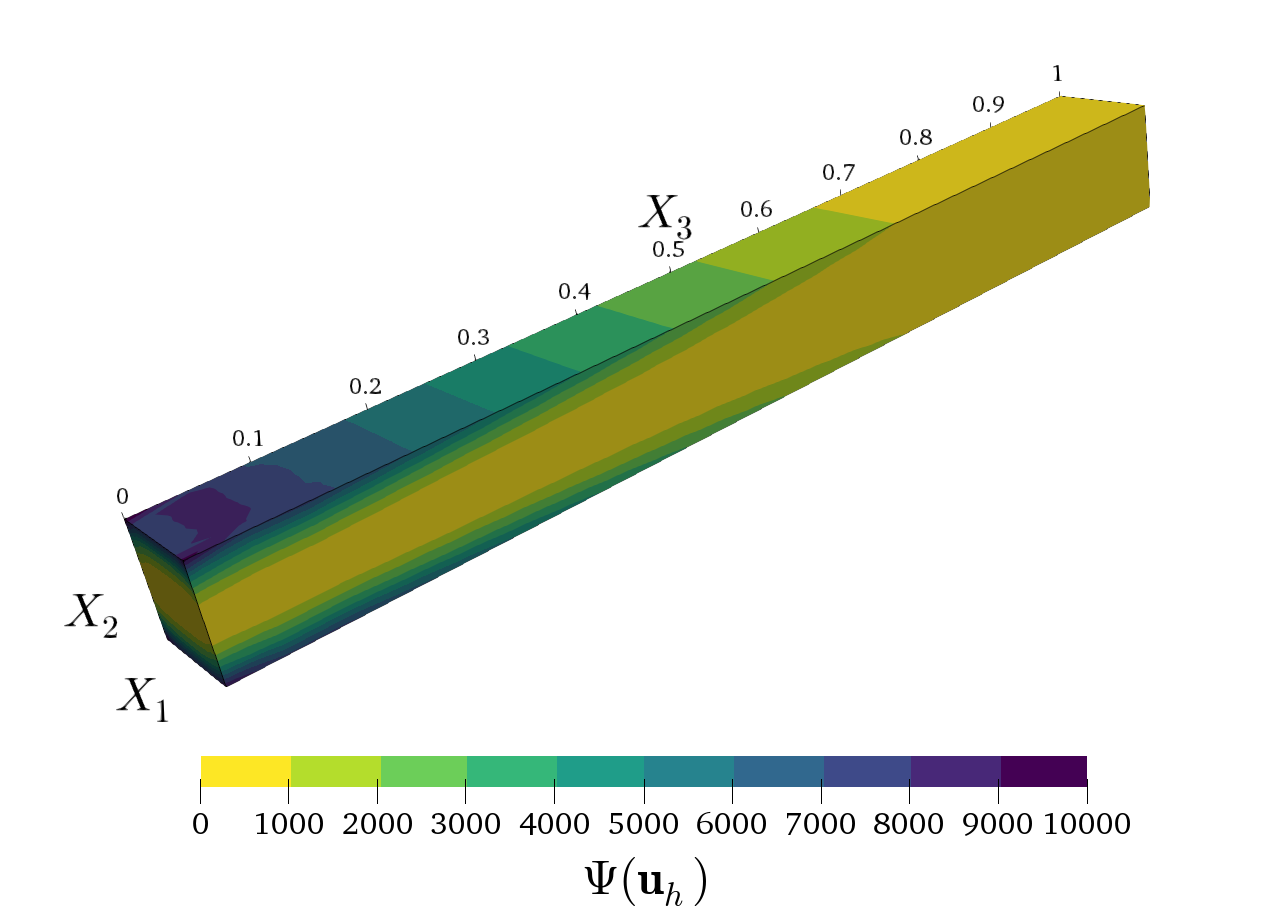}\\
  (b)
\end{subfigure}
\caption{Example 5a: Mesh refinement study of the cantilever simulation
outlined in Section~\ref{sec:hyperelasticity}: (a) Convergence of the relative
error in the functionals of the DGFEM approximation; (b) The computed internal
strain energy.}
\label{fig:cantilever}
\end{figure}

\begin{figure}[t!]
\centering
\begin{subfigure}{.5\textwidth}
  \centering
  \includegraphics[width=1.\linewidth]{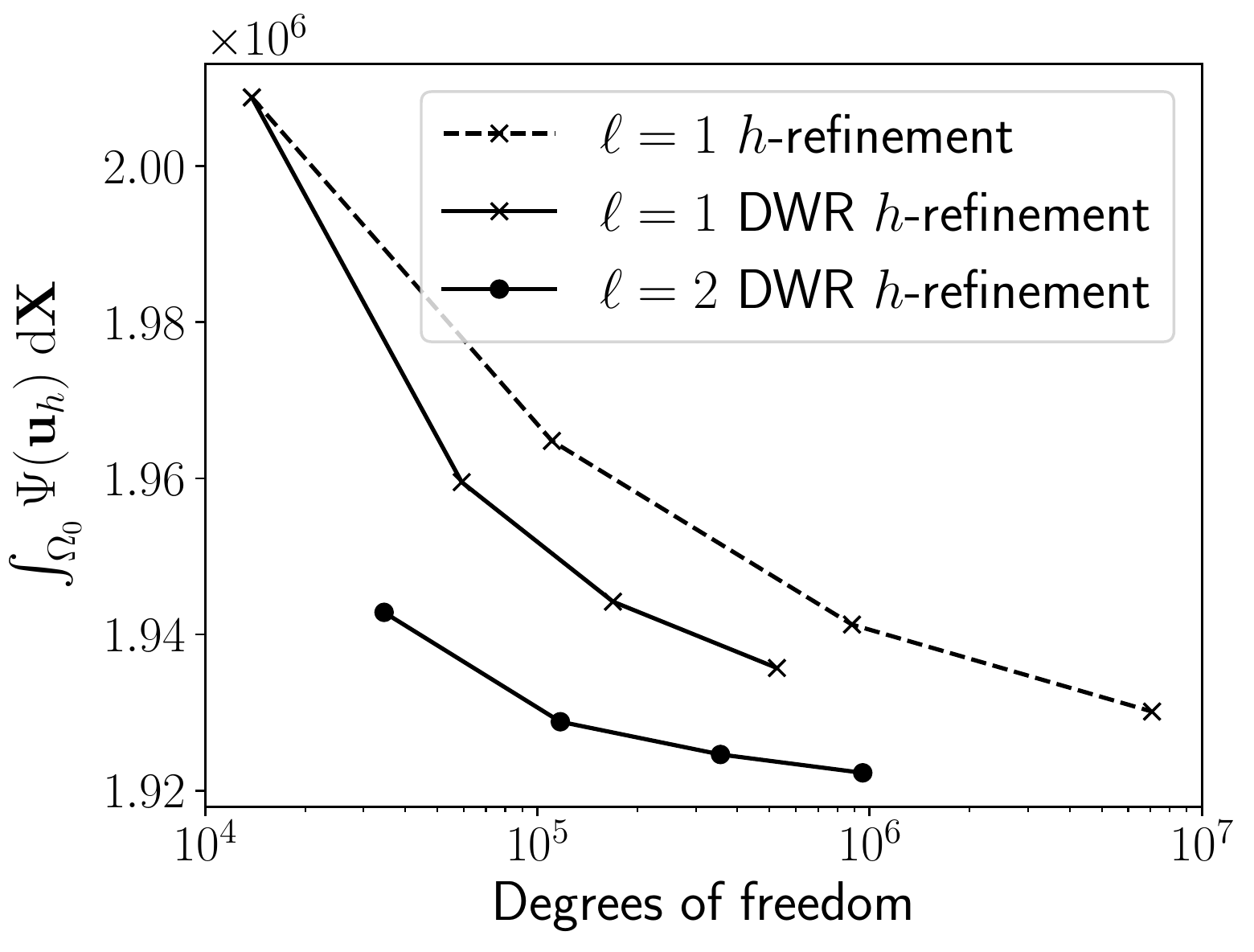}\\
  (a)
\end{subfigure}
\begin{subfigure}{.3\textwidth}
  \centering
  \includegraphics[width=1.\linewidth]{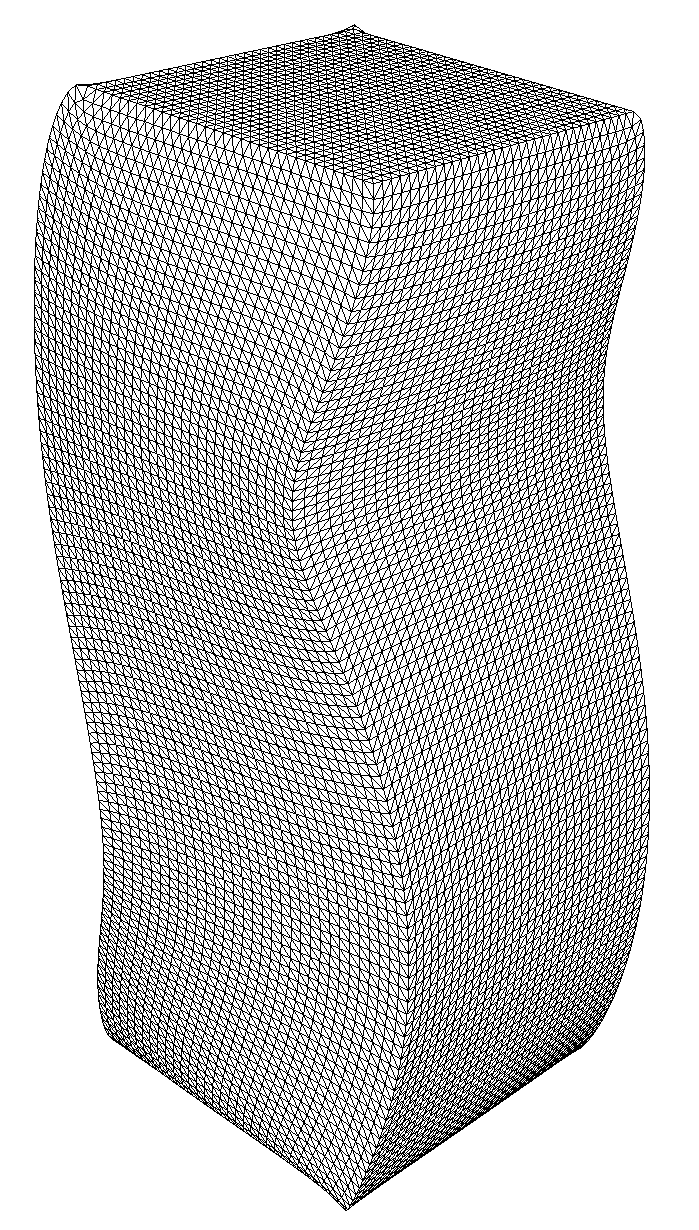}\\
  (b)
\end{subfigure}%
\begin{subfigure}{.3\textwidth}
  \centering
  \includegraphics[width=1.\linewidth]{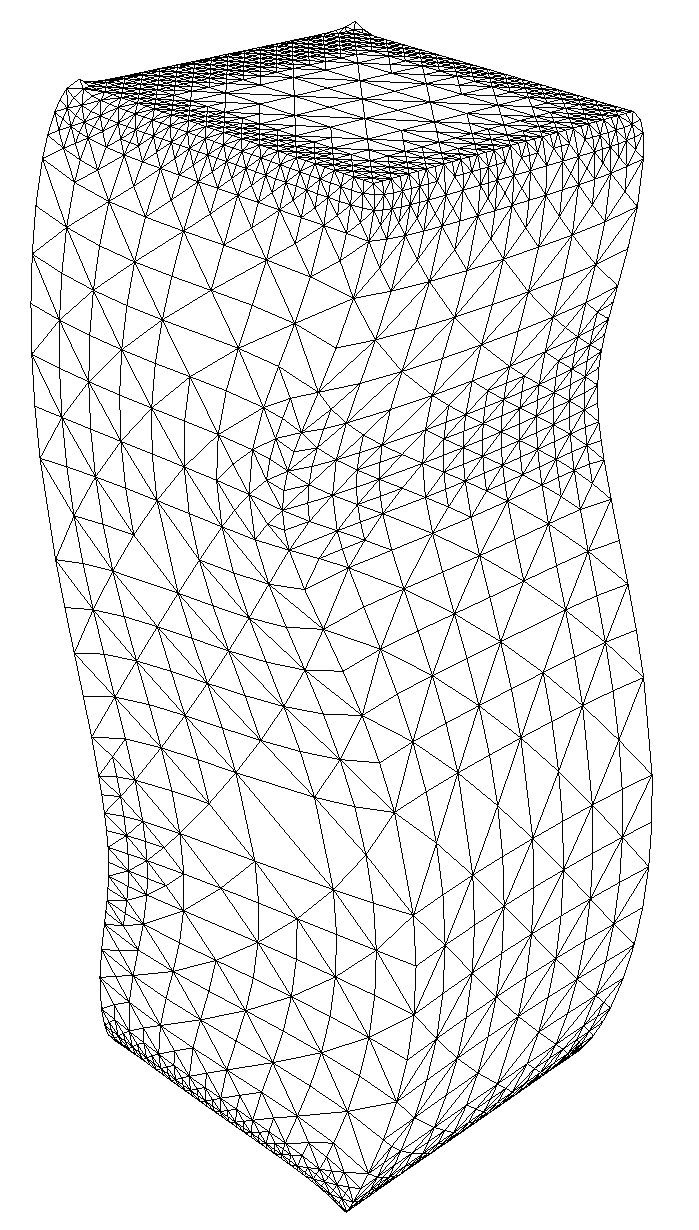}\\
  (c)
\end{subfigure}
\caption{Example 5b: Mesh refinement study of the compressible
  neo-Hookean elasticity model: (a) Convergence of the quantity
  of interest, $\int_{\Omega_0} \Psi(\mathbf{u}_h) \; \mathrm{d}\mathbf{X}$;
   (b) Deformed configuration generated by employing $\ell = 1$ and uniform $h$--refinement
  comprising~\num{7077888} degrees of freedom; (c) Deformed configuration
  generated using adaptive DWR $h$--refinement with $\ell = 2$ giving rise to \num{953220}~degrees of
  freedom.}
\label{fig:hyperelastic_buckling}
\end{figure}

Clearly we observe that equation~\eqref{eq:hyperelasticity} can be rewritten with
$\mathcal{F}^v(\mathbf{u}, \nabla_{\mathbf{X}} \mathbf{u}) \equiv \mathbf{P}$, cf.
Listing~\ref{code:ufl_neohookean}; thereby,
the DGFEM framework outlined previously may be directly applied. 
We highlight that the flexibility of this approach naturally allows us to consider other
potential strain energy models, such as Saint Venant--Kirchhoff, Ogden, and 
Mooney--Rivlin~\cite{holzapfel2000nonlinear}.

To test the DGFEM solver we repeat the numerical example~\num{4.1} presented
in~\cite{Noels2006}. To this end, we simulate the small strain deformation of
the beam $\Omega_0 = (0, \beta)^2 \times (0, L)$, where $\beta = 0.1$, $L=1$, 
with the
material parameters $E = \num{200e9}$ and $\nu = 0.3$. There is no applied 
body force, i.e., $\mathbf{B} = \mathbf{0}$, the beam is clamped on
the near face, $\partial\Omega_{0,D} = (0, \beta)^2 \times \{0\}$,
$\mathbf{u}|_{\partial\Omega_{0,D}} = \mathbf{0}$, an external distributed
force is applied to the far face
\begin{equation*}
\mathbf{T} = 
\begin{cases} \left(0, 10^4/\beta^2, 0 \right)^\top &\mbox{if } X_3 = L, \\ 
\mathbf{0} &\mbox{otherwise},
\end{cases} 
\end{equation*}
and $\partial\Omega_{0,N} = \partial\Omega_{0} \setminus
\partial\Omega_{0,D}$. The system has a known analytical solution for 
the tip deflection
$\mathbf{u}_2(\sfrac{\beta}{2}, \sfrac{\beta}{2}, L) = \num{2e-3}$ and 
the internal strain energy
$\int_{\Omega_0} \Psi(\mathbf{u}) \; \mathrm{d}\mathbf{X} = 10$. Convergence
of the relative errors of the DGFEM approximation of this problem are shown
in Figure~\ref{fig:cantilever}. The code to generate these results is available 
in the file~\url{dg_cantilever.py}.

\subsection{Example 5b: Near--Incompressible Hyperelastic Buckling}
\label{sec:hyperelasticity2}
To demonstrate the capability of the DGFEM framework applied to larger three--dimensional
problems we draw inspiration from a numerical example of an
elastic body buckled in a compressive state, cf. Section~6.7
in~\cite{Eyck2006}. Thereby, we set $\Omega_0 = (0,\sfrac{1}{2})^2 \times (0,
2)$, $E = \num{1e8}$, $\nu = 0.46$, $\mathbf{B} = \mathbf{T} = \mathbf{0}$ and
prescribe $\mathbf{u}|_{X_3 = 0} = \mathbf{0}$ and $\mathbf{u}|_{X_3=2} = (0,
0.07, -0.5)^\top$. Using the proposed DGFEM framework the
resulting deformed mesh configuration is shown in
Figure~\ref{fig:hyperelastic_buckling}. Here, we have employed both uniform
mesh refinement, as well as adaptive refinement of the computational
mesh based on employing a DWR {\em a posteriori} error indicator, where
we select the internal strain energy $\int_{\Omega_0}
\Psi(\mathbf{u}) \; \mathrm{d}\mathbf{X}$ to be the quantity of interest. 
As expected, the adaptive refinement algorithm selects regions to enrich the
computational mesh near the edges of
the Dirichlet boundary, where large changes in the stress occur, as well as
regions under compression, which lead to large changes in the internal strain energy.
The implementation of this example is provided in the
file \url{dg_compression.py}.
%





\section{Concluding remarks} \label{sec:con_rem}

In this article we have exploited the use of symbolic algebra for the automatic
computation of DGFEMs for the numerical approximation of general systems of
nonlinear PDEs. In particular, we have proposed and implemented a class
structure in order to allow for the specification of a given DGFEM in a clear
and concise manner. While the examples we have presented have primarily focused
on flow problems, we stress that the generality of this approach allows for the
treatment of a wide range of PDEs stemming from diverse application areas.
Indeed, in the final two examples, we have included the DGFEM approximation of the
indefinite Maxwell problem and a hyperelasticity problem, respectively. The exploitation of
the software developed within this article allows the user to build up the DGFEM
discretisation of systems of multi-physics PDEs in a simple and concise manner;
as an example, we have employed the software here for the DGFEM approximation of
microwave power assisted chemical vapour deposition reactor employed for the
manufacture of synthetic diamond; see~\cite{cvd_paper,nate_phd}. Moreover, the
code can easily be extended to other problems, and indeed other DGFEM schemes,
so that users can tailor the software to their own applications.

\bibliographystyle{amsplain}
\bibliography{literature}
\end{document}